\documentclass[11pt]{article}
\usepackage{amsmath,amsfonts,latexsym,amssymb,amscd, graphicx,pictexwd, mathrsfs, dsfont, color,mathtools,hyperref,subfig,caption,wasysym}
\usepackage[title]{appendix}

\newcommand{\w}{{\mathrm w}}
\newcommand{\y}{{\mathrm y}}

\newcommand{\Exp}{{\rm Exp}}
\newcommand{\Poi}{{\rm Poi}}

\renewcommand{\P}{{\mathbb P}}
\newcommand{\E}{{\mathbb E}}

\newcommand{\R}{{\mathbb R}}
\newcommand{\bN}{{\mathbb N}}
\newcommand{\D}{{\mathbb D}}
\newcommand{\bZ}{{\mathbb Z}}

\newcommand{\ZZ}{{\mathbb Z}}

\newcommand{\PP}{{\mathbf P}}

\newcommand{\cc}{{\mathbf c}}
\newcommand{\zz}{{\mathbf z}}
\newcommand{\kk}{{\mathbf k}}

\newcommand{\xx}{{\mathbf x}}

\newcommand{\yy}{{\mathbf y}}

\newcommand{\cF}{{\mathcal F}}
\newcommand{\cS}{{\mathcal S}}
\newcommand{\cM}{{\mathcal M}}

\newcommand{\cI}{{\mathcal I}}

\newcommand{\cD}{{\mathcal D}}

\renewcommand{\H}{{\mathit H}}

\newcommand{\1}{{\mathds{1}}}

\renewcommand{\d}{{\rm d}}
\newcommand{\per}{{\rm per}}

\newcommand{\dst}{\stackrel{\mbox{\tiny$dist.$}}{=}}
\newcommand{\df}{\stackrel{\mbox{\tiny$def.$}}{=}}
\newcommand{\as}{\stackrel{\mbox{\tiny$a.s.$}}{=}}

\newtheorem{thm}{Theorem}%[section]

\newtheorem{lem}{Lemma}

\newtheorem{Def}{Definition}

\newtheorem{hyp}{Hypothesis}
\numberwithin{equation}{section}
\usepackage{comment}
\begin{document}
\title{Sprinkled Decoupling for Hammersley's Process}  

\author{Leandro P. R. Pimentel and Roberto Viveros}
\date{\today}
\maketitle

\begin{abstract}  
In this article we prove a sprinkled decoupling inequality for the stationary Hammersley's interacting particle process \cite{AlDi,Ha}. Inspired by \cite{BaTe,HiKiTe}, we apply this inequality to study two distinct problems on the top of this particle process. First, we analyze a detection problem, demonstrating that a fugitive can evade  particles, provided that their jump range is sufficiently large. Second, we show that a random walk in a dynamic random environment exhibits ballistic behavior with respect to the characteristic speed of the particle system, under a weak assumption on the probability of being away of this critical speed.

\end{abstract}

\section{Models and Results}
\subsection{Introduction}
Sprinkling is a technique used to control the decay of correlations (decoupling) by introducing small perturbations. In summary, instead of considering two fixed increasing events \( A \) and \( B \), we introduce slightly ``sprinkled'' versions of these events, denoted \( A_\epsilon \) and \( B_\epsilon \), which occur with slightly higher probabilities due to the added randomness, and obtain a decoupling inequality of the form:
\[
\mathbb{P}(A \cap B) \leq \mathbb{P}(A_\epsilon ) \mathbb{P}(B_\epsilon) + \text{error}.
\]
The addition of the sprinkling parameter \( \epsilon \) accounts for the extra randomness introduced to ensure better connectivity or event occurrence. This formulation is particularly useful in percolation theory of correlated systems, where small increments in probability can lead to large-scale connectivity changes. See for instance the pioneer works \cite{PoTe,SiSz,Sz}, where sprinkling is used to control dependencies in a model of random interlacements. For more recent examples of applications of the sprinkling technique for percolation theory on Poisson cylinders and Gaussian fields we refer to \cite{AlTe,Mu}.

On the subject of conservative particle process, such as independent symmetric simple random walks \cite{BlHiSaSiTe,HiHoSiSoTe} and symmetric simple exclusion process \cite{BaTe,HiKiTe}, sprinkled decoupling has been applied to study detection problems and random walks on dynamic random environment. In both cases we have positive association and a slow decay of correlations, which makes the use of sprinkling a key tool to apply multiscale renormalization schemes. For the symmetric simple exclusion process, the symmetry of the dynamics played a crucial role in the proof methods (interchange process). The novelty of this article lies in establishing a sprinkled decoupling inequality (Theorem \ref{Sprinkling}) for a stationary, conservative and totally asymmetric particle process, known shortly as Hammersley’s process. This is achieved through a new decoupling approach that is based on geometric properties of a companion last-passage percolation model.

The description of our underlying particle process is quite simple: particles are initially distributed as a homogeneous Poisson point process of rate $\lambda$ on $\R$; then a particle at $z\in\R$ jumps to the left at rate given by the distance to its closest particle to the left, say $z^\ast$, and the new location is chosen uniformly between $z$ and $z^\ast$ (see Figure 1). Particles interact because they cannot pass through each other, and the distribution of  particle locations stay invariant for all times (stationary regime). The hydrodynamic limit $u(x,t)$ of the density profile of Hammersley's process is described by a hyperbolic conservative law, also known as Burgers equation: $\partial_t u+\partial_xf(u)=0$, where $f(u)=-1/u$ \footnote{For the heuristics behind, see the derivation of equation (9) in \cite{AlDi} for the integrated density $U=\int u dx$.}. Equations of this type often describe transport phenomena. As $\partial_x f(u)=f'(u)\partial_x u$, the directional derivative of $u$ in the direction $(1,f'(u))$ vanishes. Hence $u$ is constant on each line of the form $x=x_0+tf'(u(x_0))$, known as the characteristics lines, and one can implicitly solve the equation whereby $u(x_0+tf'(u(x_0)),t)=u(x_0)$. In fact, this is the classical method of the characteristics to solve a partial differential equation. If $u(x,0)=\lambda$ we get the slope $\rho=f'(\lambda)=\lambda^{-2}$, which will play a fundamental rule in our geometric route to decoupling. An important feature is that a disturbance  made in the initial data travels along the characteristic of the equation. In terms of the interacting particle process, this means that  space-time correlations also travel along the characteristics \footnote{See for instance equation (2.4) in the proof of Theorem 2.1 of \cite{CaGr2}.}.          
\begin{figure}\label{FigHammersley}
\begin{center}
\includegraphics[width=10cm, height=6cm]{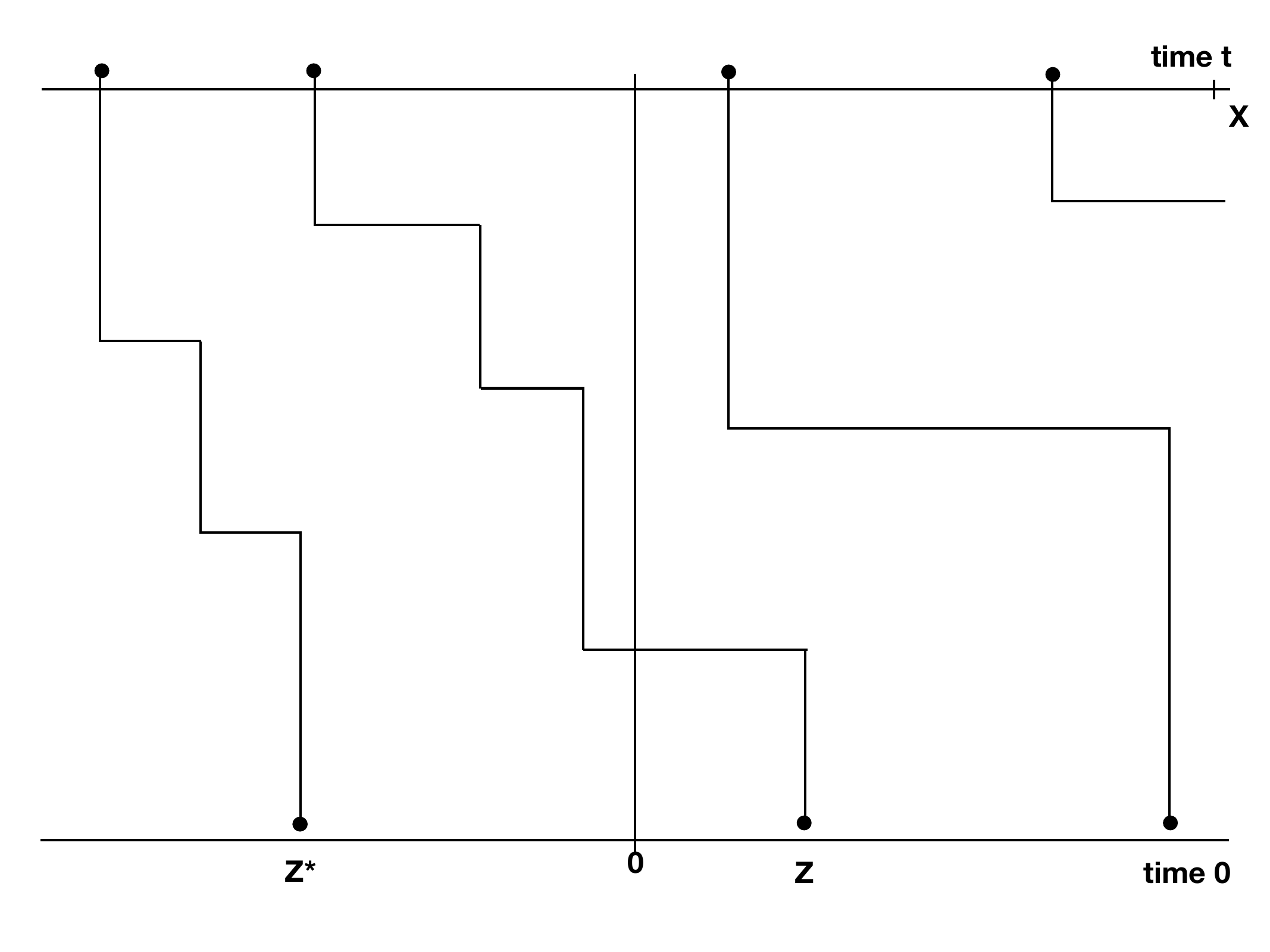}
\caption{Particles are represented by bullets ($\bullet$). Time goes up and trajectories according to Hammersley's process are represented as broken lines attached to each particle.}
\end{center}
\end{figure}

The sprinkled decoupling inequality that we obtain for Hammersley's process has a difference compare to the one obtained by Baldasso and Texeira \cite{BaTe}, for the symmetric  simple exclusion process. In their context, the inequality involves monotone functions with support on space-time boxes whose distance is large compare to the  perimeters of the boxes. In ours, we require that the distance is large compare to the perimeters multiplied by the reciprocal of the sprinkling parameter. Nevertheless, we show that this stronger requirement is not a big issue to adapt the multiscale renormalization schemes developed in \cite{BaTe} and \cite{HiKiTe}. As applications, this allows us to study two interesting problems on the top of Hammersley's process, which we describe in the next two paragraphs.

\paragraph{\bf Detection Model} Suppose that  particles are seen as detectors that are capable of detecting all targets within some distance $r>0$ from their location. Targets can  jump up to a certain range $N\geq 1$ only at discrete times, and can predict the future movement of all detectors. Since particles eventually get  arbitrarily close,  it is not hard to see that, if the target jump range $N$ is not large compared to the detection range $r$, then a target eventually gets detected with probability one. We will prove that, for sufficiently large $N$, a target starting at the origin can avoid detection forever with positive probability (Theorem \ref{Detection}).

\paragraph{\bf Random Walk in Dynamic Random Environment}
Fix two parameters $p_\bullet,p_\circ\in(0,1)$. The random walk starts at the origin and jumps at discrete time. At the moment of the jump, it checks if there is a particle within distance $1/2$ from it. If it does, then it chooses to jump to the right with probability $p_\bullet$, or to the left with probability $1-p_\bullet$. If not, then it chooses to jump to the right with probability $p_\circ$, or to the left with probability $1-p_\circ$. We will prove that the random walk exhibits ballistic behavior with respect to the velocity $\rho=\lambda^{-2}$ , under a weak assumption on the probability of being away of this critical value (Theorem \ref{Ballistic}). For rigorous results regarding the long-time behavior of random walks in symmetric dynamic random environment, such as independent random walks or the symmetric exclusion process, we address the reader to \cite{HiHoSiSoTe,HiKiTe,HuSi1,HuSi2}. For asymmetric environments, such as Hammersley's process or the totally asymmetric simple exclusion process, it is still  conjectured the existence of a velocity, and Gaussian fluctuations away from the characteristic line. Along this critical line, it is expected the emergence of super-diffusive fluctuations lying within the Kardar-Parisi-Zhang universality class. This behavior was first predicted in the physics literature by T. Bohr and A. Pikovsky \cite{BoPi}, where numerical simulations and heuristics arguments pointed it out. The decoupling inequality is the first step of an ongoing project of the authors, where the focus is the long-time behavior of random walks in asymmetric dynamic environments. 
          
\subsection{Hammersley's Process (HP)}
Recall that we start the stationary Hammersley's process (HP) with a configuration of particles distributed as a Poisson point process on $\R$ with parameter $\lambda>0$, and that a particle at $z\in\R$ jumps to the left at rate given by the distance to its closest particle to the left, say $z^\ast$, and the new location is chosen uniformly between $z$ and $z^\ast$ \cite{AlDi,Ha}. We represent a particle configuration as a locally finite counting measure $\eta_t^\lambda$ on $\R$, so that $\eta^\lambda_t(A)$ gives the number of particles in $A\subseteq\R$ at time $t\geq 0$. Then  $(\eta^\lambda_t)_{t\geq 0}$ is a time stationary Markov process \cite{AlDi,CaPi2}: $\eta^\lambda_t\dst \eta^\lambda_0$ for all $t\geq 0$. The stationary HP can be constructed so that time is running from $-\infty$ to $+\infty$. Let $\cM$ denote the space of locally finite counting measures on $\R$ and $\cD$ the space of cadlag trajectories $\eta:\R\to \cM$. We have a partial order $\eta_t\leq\zeta_t$ defined on $\cM$ as $\eta_t(A)\leq \zeta_t(A)$ for all $A\subseteq\R$, which induces a partial order defined on $\cD$ by $\eta\leq \zeta$ if $\eta_t\leq\zeta_t$ for all $t\in \R$. This particle process is attractive in the sense that if $\eta_s\leq\zeta_s$ for some $s\in\R$, and we run the particle processes with the same clocks (basic coupling), then $\eta_t\leq \zeta_t$ for all $t\geq s$. For more details on basic properties of Hammersley's process  we address the reader to the introductory text \cite{CaPi2}. A function $f:\cD\to\R$ is non-decreasing (non-increasing) if $\eta\leq\zeta$ implies that $f(\eta)\leq f(\zeta)$ ($f(\eta)\geq f(\zeta)$). We say that $f:\cD\to\R$ has support in a space-time box $B=[x,y]\times[s,t]\subseteq \R^2$ if for every $\eta,\zeta\in\cD$,
$$\eta_r(A)=\zeta_r(A)\,\,\forall\,\,A\subseteq [x,y]\,,\,r\in[s,t]\,\,\implies\,\,f(\eta)=f(\zeta)\,.$$   
Let \( \mathbb{P}^\lambda \) denote the probability measure on the space of càdlàg trajectories \( \eta :\R \to \mathcal{M} \), induced by the stationary HP with parameter \( \lambda > 0 \). Let \( \mathbb{E}^\lambda \) be the associated expectation operator. 

For $B=[x,y]\times[s,t]$ we set $\per(B)\df y-x+t-s$, and let  
$$\d(B_1,B_2)\df d_h+d_v\mbox{ where }\begin{cases}d_h\df\inf\{|x_1-x_2|\,:\,(x_i,t_i)\in B_i\,,i=1,2\} \\
d_v\df\inf\{|t_1-t_2|\,:\,(x_i,t_i)\in B_i\,,i=1,2\}\,.\end{cases}$$ 
Thus $d_h$ and $d_v$ are called the horizontal and vertical distances, respectively, between $B_1$ and $B_2$. In Theorem \ref{Sprinkling} below we consider functions $f_1$ and $f_2$ with support on a space-time box $B_1$  and $B_2$, respectively.

\begin{thm}\label{Sprinkling}
Let $\lambda,\lambda'\in I=[\lambda_1,\lambda_2]\subseteq(0,\infty)$ and $\epsilon=\left|\lambda^{-2}-\lambda'^{-2}\right|$. There exist $\epsilon_I,c_I,C_I>0$, that depends on $I$, such that for all  $\epsilon\in(0,\epsilon_I)$ and 
$$\d(B_1 ,B_2)>C_I\epsilon^{-1}\big(\per(B_1)+\per(B_2)\big)\,,$$
if $\lambda<\lambda'$ and $f_i$ are non-decreasing then
$$\E^\lambda\left[f_1f_2\right]\leq\E^{\lambda}\left[f_1\right]\E^{\lambda'}\left[f_2\right]+10 e^{-c_I\epsilon^4 \d(B_1 ,B_2)}\,,$$
if $\lambda'<\lambda$ and $f_i$ are non-increasing then 
$$\E^\lambda\left[f_1f_2\right]\leq\E^{\lambda}\left[f_1\right]\E^{\lambda'}\left[f_2\right]+10 e^{-c_I\epsilon^4 \d(B_1 ,B_2)}\,.$$ 
\end{thm}

The difference with Theorem 1.5 of Baldasso and Texeira \cite{BaTe}, in the context of the simple symmetric exclusion process, is that we require the distance between boxes to be as large as the sum of the respective perimeters times the reciprocal of the sprinkling parameter $\epsilon$. To the best of our knowledge, Theorem \ref{Sprinkling} is the first sprinkled decoupling inequality for a totally asymmetric particle system.

\subsection{Detection by Hammersley's Process}

Suppose that in \( \mathbb{R}^d \) a moving object (e.g., an intruder) attempts to traverse without being detected by a detection system \( \Phi \) which consists of randomly placed particles, each with a detection radius. Each point \( x \in \Phi \) represents a sensor or obstacle and has a detection range modeled as a ball of radius \( r > 0 \), forming the detected region:

\[
D = \bigcup_{x \in \Phi} B(x, r),
\]
where \( B(x, r) \) is the closed Euclidean ball of radius \( r \) centered at \( x \). The object follows a path \( \gamma: [0, \infty) \to \mathbb{R}^d \) and the goal is to find an infinite trajectory \( \gamma \) that avoids the detected region \( D \), meaning: \[
\gamma(t) \notin D, \quad \forall t \geq 0.
\]
If such a path exists, we say the object escapes detection. This problem automatically translates into the percolation of free space \( \mathbb{R}^d \setminus D \). The study of percolation in the vacant set has been a significant topic in continuum percolation theory. This research primarily focuses on understanding the conditions under which the unoccupied regions (vacant set) form an infinite connected component, allowing an object to traverse without intersecting any obstacles. For instance, see \cite{AhTaTe} and \cite{TyWi} for the study of the vacant set when $\Phi$ is given by the realization of a Poisson point process.
\begin{figure}\label{FigDetection}
\begin{center}
\includegraphics[width=10cm, height=6cm]{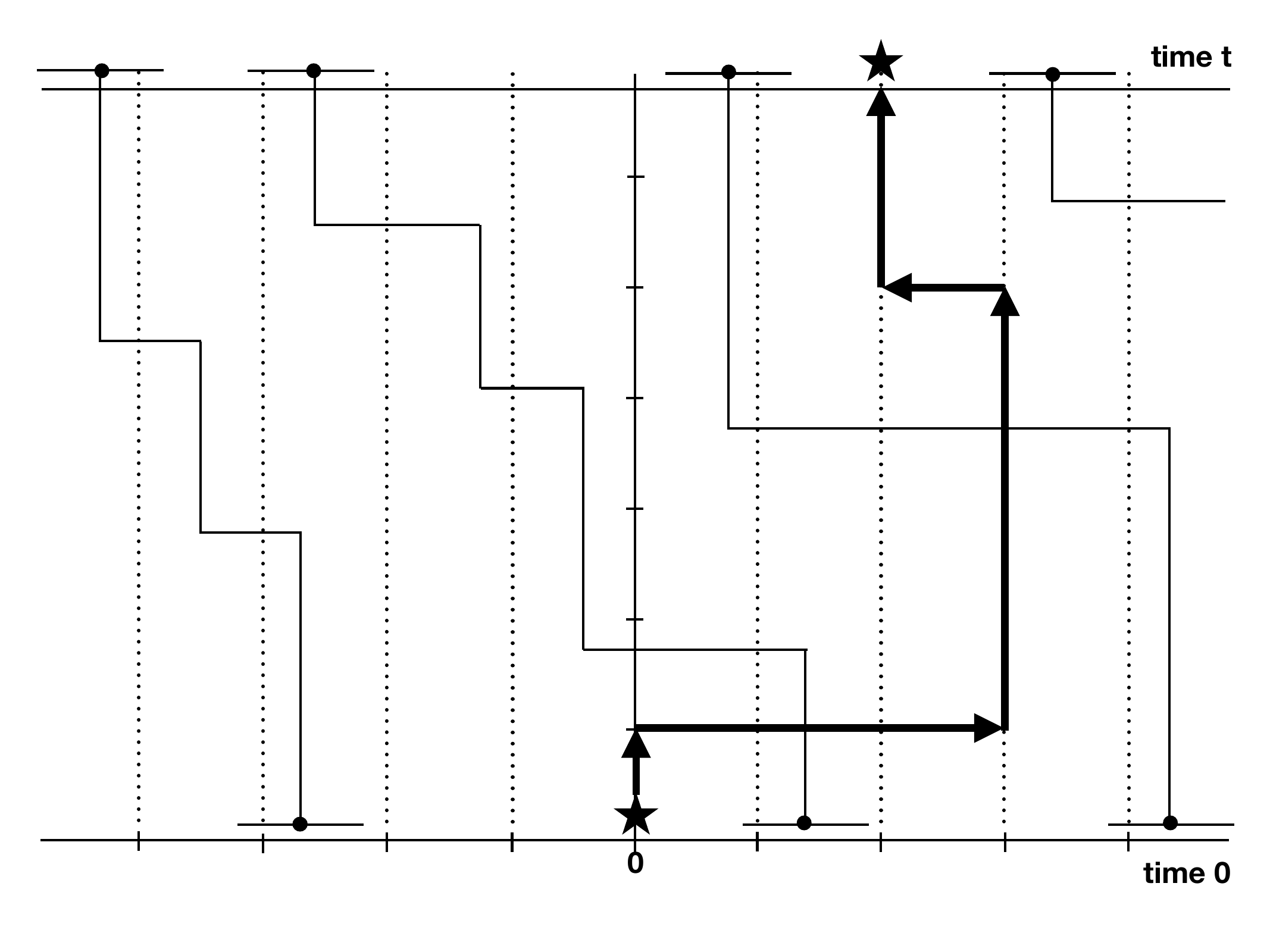}
\caption{Line segments centered at $\bullet$ represent detectors. The clairvoyant target starting at the origin is represented by $\bigstar$. A possible escape route up to time $t$ (with $N\geq 3$) is drawn with arrows.}
\end{center}
\end{figure}

In this paper, as an application of Theorem \ref{Sprinkling}, we consider the following detection problem. Suppose we have a collection of detectors that are moving over the one-dimensional real line according to Hammersley's process, and that detectors detect everything that is within distance $r>0$ from it. A target can make jumps of size $\{-N\Delta_s,\dots,N\Delta_s\}$, where $\Delta_s>0$ is a real space parameter, and $N\in\{1,2,\dots\}$ is fixed,  at times $t_n=n\Delta_t$ for $n\geq 1$, where $\Delta_t>0$ is a real time parameter. It also assumed that the target knows the future movement of all detectors, and the main question is whether it can escape detection forever with positive probability (Figure 2). It is not hard to see that, with probability one, eventually consecutive detectors gets arbitrarily close to each other. Hence, for $\Delta_s=1$, $r\geq 1/2$ and $N= 1$, with probability one, a target starting at the origin eventually gets detected. Here we prove that, for fixed $\lambda,r,\Delta_s,\Delta_t>0$, a phase transition occurs in the probability of detection as $N$ varies.
               
\begin{thm}\label{Detection}
There exists $N_0$ such that if $N\geq N_0$, then the probability that a target starting at the origin never gets detected by some detector is positive. 
\end{thm}

\subsection{Random Walk on Top of the Hammersley's Process}

Random walks in random environments (RWRE) refer to stochastic processes where a walker moves in a space where the transition probabilities are random variables. More specifically
\[
\mathbb{P}_\omega(X_{n+1} = x+e_i | X_n = x) = \omega(x,n, e_i),
\]
where $\omega$ is an underlying random environment drawn from a probability space $\Omega$. This environment is usually fixed before the walk begins (quenched randomness), and the idea is to analyze the conditional behavior of the walker given the specific realization of $\omega$. Random walks in random environments have been widely studied, in both static (when doesn't depend on $n$) and dynamic settings. A particularly interesting class of models arises when the environment is given by a conservative particle system, such as independent random walks \cite{HiHoSiSoTe} or the exclusion process \cite{HiKiTe}. In this case, the correlations in the environment are strong and the use of sprinkling, which helps in decoupling space-time events, has been crucial in analyzing these systems.

In this paper, also as an application of Theorem 1, we follow the multiscale strategy developed in \cite{HiHoSiSoTe, HiKiTe} for dynamic random environments, and use our inequality to implement their framework for a random walk that evolves on top of  Hammersley’s process. More specifically, we prove that under a suitable condition, the random walk remains with high probability on one side of a straight line aligned with the characteristic direction of the environment. Depending on the parameters, this side may be to the left or right of the characteristic. While this does not establish ballistic behavior in the strong sense (such as a law of large numbers), it does imply that the walker makes linear progress in a given direction.

Let $I_x\df(x-1/2,x+1/2)$ denote the open interval of diameter one around $x$, and fix $p_\bullet\,,\,p_\circ\in(0,1)$. Define the random walk $(X_n)_{n\geq 0}$ on $\bZ$ starting at $X_0=0$ and evolving on the top of $\eta\in\cD$ as follows:  
\begin{itemize}
\item If $\eta_n(I_{X_n})>0$ then $\begin{cases}X_{n+1}=X_n +1\mbox{ with probability }p_\bullet\\  X_{n+1}=X_n -1\mbox{ with probability }1-p_\bullet\end{cases}$;
\medskip

\item If $\eta_n(I_{X_n})=0$ then $\begin{cases}X_{n+1}=X_n +1\mbox{ with probability }p_\circ\\ X_{n+1}=X_n -1\mbox{ with probability }1-p_\circ\end{cases}$.

\end{itemize}
Let $P^\eta$ denote the probability law of $(X_n)_{n\geq 0}$, given $\eta\in\cD$, also known as the quenched law of the random walk, and let 
$$\PP^\lambda(\cdot)\df\E^\lambda\left[P^\eta(\cdot)\right]\,,$$
also known as the annealed law of the random walk.

 We consider the coupled family \( \{(X_t^\w)_{t \geq 0},\; \w \in \mathbb{R}^2\} \) constructed as in Section 3.1 of \cite{HiKiTe}. The position of the random walk that starts at $\w = (x,n) \in 
\mathcal{L} = (2\mathbb{Z})^2 \cup \left( (1,1) + (2\mathbb{Z})^2 \right)$ after one step is
\[
X_1^\w =
\begin{cases}
x + 2\,\mathbf{1}_{\{ U_\w \leq p_\bullet \}} - 1, & \text{if } \eta^\lambda_n(I_x) > 0, \\
x + 2\,\mathbf{1}_{\{ U_\w \leq p_\circ \}} - 1, & \text{if } \eta^\lambda_n(I_x) = 0,
\end{cases}
\]
where $\{U_\w: \w \in \mathcal{L}\}$ is an i.i.d. collection of uniform random variables in $[0,1]$. For any integer \( m \geq 1 \),  the walk is defined by induction:
\[
X_m^\w = X_1^{(X_{m-1}^\w, m - 1)},
\]
which recursively constructs the trajectory starting from \( \w \in \mathcal{L} \) and then is extended to all possible starting points in \( \mathbb{R}^2 \) assigning each point to the closest previously defined trajectory located vertically above it (see the details in Section 3.1 of \cite{HiKiTe}). This collection of random walks satisfies two fundamental properties:
\begin{itemize}
    \item \textbf{Coalescence:} if two random walks meet at some vertex, they will follow the same trajectory from that point onward;
    \item \textbf{Monotonicity:} if one walk starts to the right of another, it will remain to the right for all future times.
\end{itemize}

We can add to this construction a coupling of HP at different intensities. More precisely, using the basic coupling of the family of measures \( \{ \mathbb{P}^\lambda : \lambda > 0 \} \), we can define a new coupled family of random walks \( (X_t^{\w, \lambda})_{t \geq 0},\; \w \in \mathbb{R}^2 \), where \( \lambda \) denotes the intensity of the underlying Poisson point process that initializes the environment. In this setting, each trajectory \( (X_t^{\w, \lambda})_{t \geq 0} \) evolves in the dynamic environment \( \eta^\lambda \) given by the stationary Hammersley's process with parameter \( \lambda>0 \). The use of the basic coupling ensures that if \( \lambda < \lambda' \), then the environment \( \eta_t^\lambda \) is contained in \( \eta_t^{\lambda'} \) for all \( t \in\R \), almost surely. As a consequence of the monotonicity of the walks with respect to the environment and the attractiveness of the basic coupling, it readily follows that if \( p_\circ < p_\bullet \), then
\[
X_t^{\w, \lambda} \leq X_t^{w, \lambda'} \quad \text{for all } t \geq 0,\quad \text{whenever } \lambda < \lambda'.
\]

For $\w=(w_1,w_2)\in\R^2$ let $\pi_1(\w)\df w_1$ and $\pi_2(\w)\df w_2$, and denote by $\D_1$ the set of all $\w\in\R^2$ such that $|\pi_1(\w)|+|\pi_2(\w)|\leq 1$. As in \cite{HiKiTe}, we consider the following events 
$$A_{H,\w}(v, \lambda)\df\left\{\exists\,\y\in\left(\w+[0,H)\times\{0\}\right)\mbox{ s.t. } X^{\y,\lambda}_H-\pi_1(\y)\geq v H \right\}\,$$
and
$$\tilde{A}_{H,\w}(v,\lambda)\df\left\{\exists\,\y\in\left(\w+[0,H)\times\{0\}\right)\mbox{ s.t. }X^{\y,\lambda}_H-\pi_1(\y)\leq v H \right\}\,.$$
In what follows, we parametrize the random walk on Hammersley's process by the characteristic direction $\rho=\lambda^{-2}$, and define 
$$p_H(v,\rho)\df \sup_{\w\in\R^2}\PP\left(A_{H,\w}(v,\lambda)\right)= \sup_{\w\in\D_1}\PP\left(A_{H,\w}(v,{\lambda})\right) $$
and 
$$\tilde{p}_H(v,\rho)\df \sup_{\w\in\R^2}\PP\left(\tilde A_{H,\w}(v,\lambda)\right) =\sup_{\w\in\D_1}\PP\left(\tilde A_{H,\w}(v,\lambda)\right)\,.$$
We note that $p_H(v,\rho)$ and $\tilde{p}_H(v,\rho)$ depend on $p_\bullet\,,\,p_\circ\in(0,1)$ as well. The following quantities are the natural candidates for the speed of the random walk:
$$v_+(\rho)\df\inf\left\{v\in\R\,:\,\liminf_{H\to\infty}p_H(v,\rho)=0 \right\}\,,$$
and 
$$v_-(\rho)\df\sup\left\{v\in\R\,:\,\liminf_{H\to\infty}\tilde{p}_H(v,\rho)=0 \right\}\,.$$

If \( p_\circ \leq p_\bullet \) and $\lambda \leq \lambda'$ then $A_{H,\w}(v, \lambda) \subset A_{H,\w}(v, \lambda') $. This implies that the function \( p_H(v, \rho) \) is non-increasing in \( \rho \), which implies that \( \rho \mapsto v_+(\rho) \) is also non-increasing. Note that when \( \rho = 0 \), the system can be interpreted as having an infinitely dense initial Poisson process. In this regime, the random walk always encounters a particle at every vertex and by the {law of large numbers}, this leads to:
\[
\frac{X_H}{H} \to 2p_\bullet - 1 \quad \PP-\text{a.s.}
\]
as \( H \to \infty \). This implies that \( v_+(0) = v_-(0) = 2p_\bullet - 1\). Similarly, for \( \rho \to \infty \), we have  
\(
\lim_{\rho \to \infty} v_+(\rho) < 2p_\circ - 1.
\)
The same reasoning applies to the function \( v_-(\rho) \). In summary, both \( v_+(\rho) \) and \( v_-(\rho) \) are non-increasing functions, and their range is included in \( [2p_\circ -1 , 2p_\bullet -1].\) On the other hand, if $p_\bullet\leq p_\circ$ then $\rho\mapsto v_{\pm}(\rho)$ is non-decreasing and their range in included in \( [2p_\bullet -1 , 2p_\circ -1].\) For $p_\circ<p_\bullet$ we have that $\{\rho>0\,:\,v_-(\rho)>\rho\}$ is non-empty and bounded from above, and $\{\rho>0\,:\,v_+(\rho)<\rho\}$ is non-empty and bounded from below. Thus we can define
$$\rho_{c-}\df\sup \{\rho>0\,:\,v_-(\rho)>\rho\}\in(0,\infty)\,,$$
and 
$$\rho_{c+}\df\inf \{\rho>0\,:\,v_+(\rho)<\rho\}\in(0,\infty)\,.$$
Since $v_{\pm}(\rho)$ is non-increasing, we also have that 
$$\rho_{c-}\df\inf \{\rho>0\,:\,v_-(\rho)<\rho\}\,,$$
and 
$$\rho_{c+}\df\sup \{\rho>0\,:\,v_+(\rho)>\rho\}\in(0,\infty)\,.$$

We adapt the multiscale analysis in \cite{HiKiTe} to fit the requirement in Theorem \ref{Sprinkling} for the decoupling inequality, and prove ballisticity of the random walk with respect to the characteristic speed, under a weak condition on the probability of being away from this critical value \footnote{When we say ballisticity with respect to some speed $v$, it means that the random walk is moving strictly faster or strictly slower than that speed almost surely. Thus, the usual definition of ballisticity corresponds to $v=0$.}.

\begin{thm}\label{Ballistic}
Assume that $p_\circ<p_\bullet$. For any $\rho<\rho_{c-}$ we have that $v_-(\rho)>\rho$ and there exist $\tilde{c}_1\in(0,\infty)$ and $\tilde{v}>\rho$ such that 
$$\tilde{p}_H(\tilde{v},\rho)\leq \tilde{c}_1\exp\left(-2\log^{3/2} H \right)\,,\,\forall\,H\geq 1\,.$$
For any $\rho>\rho_{c+}$ we have that $v_+(\rho)<\rho$ and there exist $c_1\in(0,\infty)$ and $v<\rho$ such that 
$$p_H(v,\rho)\leq c_1\exp\left(-2\log^{3/2} H \right)\,,\,\forall\,H\geq 1\,.$$
\end{thm}

\subsection{Sprinkled Decoupling for Exponential Last-Passage Percolation}
As mentioned before, the approach to decoupling is based on geometric properties of a Poisson last-passage percolation (LPP) model that is related to Hammersley's process. It turns out that the same approach can be applied to a lattice LPP model with exponential weights to prove a similar sprinkled decoupling inequality. To introduce the model denote $\ZZ_+\df\{0,1,2,\dots\}$ and assume that $\{w_{i,j}\,:\,(i,j)\in\ZZ_+^2\}$ is an array of nonnegative numbers such that $w_{0,0}=0$. Let $\Pi_{ij}$ be the set of oriented paths
$$\pi=\left((0,0)=(k_0,l_0),(k_1,l_1),\dots,(k_{i+j},l_{i+j})=(i,j)\right)\,,$$
with up-right steps: $(k_{n+1},l_{n+1})-(k_{n},l_{n})\in\{(1,0)\,,\,(0,1)\}$. Whenever we denote $(k,l)\in\pi$, it means that $(k,l)\in\{(k_0,l_0),\dots,(k_{i+j},l_{i+j})\}$. The last-passage percolation time from $(0,0)$ to $(i,j)$ is defined as 
$$(i,j)\in\ZZ_+^2\,\mapsto\,G(i,j)\df \max_{\pi\in\Pi_{ij}}\sum_{(k,l)\in\pi} w_{kl}\,,$$
and we say that $\pi\in\Pi_{ij}$ is a maximal path if $G(i,j)=\sum_{(k,l)\in\pi} w_{kl}$. To introduce randomness, we pick  $\alpha\in(0,1)$ and take the random variables $\omega_{ij}$ for $(i,j)\in\ZZ_+^2\setminus\{(0,0)\}$ independent with marginal distributions described below:
\begin{eqnarray*}   
w_{i0}&\dst&\Exp(1-\alpha)\,,\,i\geq 1\,,\\
w_{0j}&\dst&\Exp(\alpha)\,,\,j\geq 1\,,\\
w_{ij}&\dst&\Exp(1)\,,\,i,j\geq 1\,.
\end{eqnarray*}
Here $w\dst\Exp(\alpha)$ means that $w$ has the exponential distribution with rate $\alpha$. We introduce the superscript $\alpha$ to indicate parameter dependence. 

This LPP model is related to the totally asymmetric simple exclusion process (TASEP), an interacting particle system $\left(\xi^\alpha_t(x)\,,\,x\in\ZZ\right)_{t\geq 0}$ where each particle waits a rate one exponential random amount of time, then attempts to jump one step to the right and the jump is allowed only if there is no particle in the target site. The evolution takes place in $\{0,1\}^{\ZZ}$, where at each time $t\geq 0$, $\xi^\alpha_t(x)=1$ indicates that there is a particle at site $x$, while $\xi^\alpha_t(x)=0$ indicates that there is a hole (no  particle) at site $x$. Initially, particles at sites $x\in\ZZ$ obey independent Bernoulli with density $\alpha\in(0,1)$ distributions. Then we have a time-stationary process: $\xi^\alpha_t\dst\xi_0$ for all $t\geq 0$.  We condition on $\xi^\alpha_0(0)=0$ and $\xi^\alpha_0(1)=1$, and the particle initially at site one is labeled $0$, and the hole initially at the origin is labeled $0$. After this, all particles are labeled with integers from right to left, and all holes from left to right.  As explained in \cite{BaCaSe}, the collection of last-passage percolation times $\{G^\alpha(i,j)\,:\,(i,j)\in\ZZ_+^2\}$ has the same distribution as the collection of times where the $j$th particle and $i$th hole exchange positions. The TASEP also has a hydrodynamic limit described by a hyperbolic conservative law with $f(u)=u(1-u)$. In the stationary regime, the characteristic velocity is given by $f'(\alpha)=1-2\alpha$. In the two-dimensional LPP model, this corresponds to the direction $\left((1-\alpha)^2,\alpha^2\right)$, and the characteristic velocity $\beta=\beta(\alpha)\df (1-\alpha)\alpha^{-2}$. For a more detailed explanation, we refer the reader to section 3.4 in \cite{BaCaSe}.

The LPP model induces a discrete time evolution $W^\alpha=(W^\alpha_n)_{n\geq 0}$ defined as 
\begin{equation}\label{ExpLPP}
W^\alpha_n((k,l])\df G^\alpha(l,n)-G^\alpha(k,n)\,,\,\mbox{ for }(k,l]\subseteq \ZZ_+\,,
\end{equation}
which is time stationary: $W^\alpha_n\dst W^\alpha_0$ for all $n\geq 0$ (Lemma 4.2 \cite{BaCaSe}). We have a partial order on the space of possible evolutions, defined as $W\leq W'$ if $W_n((k,l])\leq W'_n((k,l])$ for all $(k,l]\subseteq \ZZ_+^2$ and $n\geq 0$. It is also true that if $w_{ij}$ and $w'_{ij}$ 
are such that $w_{00}=\tilde w'_{00}=0$, $w_{j0}\geq w'_{0j}$, $w_{i0}\leq w'_{i0}$ and $w_{ij}= w'_{ij}$ for all $i,j\geq 1$ then $W\leq W'$ (Lemma 4.3 \cite{BaCaSe}). This implies that one can couple $(W^\alpha,W^{\alpha'})$ with $\alpha'\leq \alpha$ such that $W^\alpha\leq W^{\alpha'}$. In what follows we consider two monotonic functions $f_1$ and $f_2$ of $W$ with support on space-time boxes $B_1$ and $B_2$, respectively.  
\begin{thm}\label{SprinklingLPP}
Let $\alpha,\alpha'\in I=[\alpha_1,\alpha_2]\subseteq(0,1)$, and let  $\epsilon=\left|\beta(\alpha)-\beta(\alpha')\right|$. There exist $\epsilon_I,c_I,C_I>0$, that depends on $I$, such that for all  $\epsilon\in(0,\epsilon_I)$ and 
$$\d(B_1 ,B_2)>C_I\epsilon^{-1}\big(\per(B_1)+\per(B_2)\big)\,,$$
if $\alpha'<\alpha$ and $f_i$ are non-decreasing then
$$\E^\alpha\left[f_1f_2\right]\leq\E^{\alpha}\left[f_1\right]\E^{\alpha'}\left[f_2\right]+10 e^{-c_I\epsilon^4 \d(B_1 ,B_2)}\,,$$
if $\alpha<\alpha'$ and $f_i$ are non-increasing then 
$$\E^\alpha\left[f_1f_2\right]\leq\E^{\alpha}\left[f_1\right]\E^{\alpha'}\left[f_2\right]+10 e^{-c_I\epsilon^4 \d(B_1 ,B_2)}\,.$$ 
\end{thm}

Theorem 4 can be used to apply multiscale analyze to a random walk on the top of a exponential LPP environment, and deduce a result similar to Theorem \ref{Ballistic}. In this model, for each $x\in\ZZ$ and $n\geq 0$ we introduce the derivative 
$$D^\alpha\left(x,n\right)\df G^\alpha\left(\frac{n+x}{2},\frac{n-x}{2}+1\right)-G^\alpha\left(\frac{n+x}{2}+1,\frac{n-x}{2}\right)\,,$$
and the random walk located at position $x\in\ZZ$ at time $n\geq 0$ inspects the LPP environment at the space-time position where it lies on: if $D^\alpha\left(x,n\right)>0$ , the random walk decides to jump to the right with probability $p_+$ or to the left with probability $1-p_+$; if $D^\alpha\left(x,n\right)<0$, the random walk decides to jump to the right with probability $p_-$ or to the left with probability $1-p_-$. This is part of a ongoing project of the authors, which also includes a more detailed study of the long time behavior of a random walk in last-passage percolation environments. 

\subsection{Overview of the Paper}
In Section 2, we present the graphical construction of Hammersley's process with sources and sinks, along with key results concerning the geometry of maximal paths, which serves as the relevant quantity for establishing domination at different densities. The section concludes with the proof of the sprinkled  decoupling inequality (Theorem \ref{Sprinkling}). The proof of the analog result in the exponential LPP model (Theorem \ref{SprinklingLPP}) follows the same lines and we only present the script to achieve it. In Section 3, we demonstrate that, for sufficiently large jumps, an intruder can successfully evade detection (Theorem \ref{Detection}). The random walk result (Theorem \ref{Ballistic}) is obtained in Section 4. Both results are obtained through a combination of a renormalization scheme argument and the sprinkled decoupling inequality. For the sake of completeness, in the appendix, we provide the proof of a basic large deviation result for the Poisson distribution, that is used in the proofs.

\section{Sprinkled Decoupling for Hammersley's Process} 

\subsection{Graphical Construction and Last-Passage Percolation} 
Hammersley's process can be constructed by considering a Poisson point process $\omega$ on $\R\times(0,\infty)$ with parameter $1$, independently of $\eta^\lambda_0$, and by setting that if $(x,t)\in\omega$ then at time $t$ the particle closest to $x$, among those to the right of $x$, moves to $x$. This graphical construction allow us to describe this particle system in terms of a Poisson last-passage percolation model as follows. Call a sequence $(x_1,t_1),\dots,(x_k,t_k)$ of planar points an  increasing path if  $(x_j,t_j)<(x_{j+1},t_{j+1})$ (i.e. $x_j<x_{j+1}$ and $t_j<t_{j+1}$)  for all $j=1,\dots,k-1$. The last-passage percolation time $L((z,s),(x,t))$ between  $(z,s)<(x,t)$ is the maximal number of Poisson points in $\omega$ among all increasing paths of Poisson points lying in the rectangle $(z, x ]\times(s, t]$. For $x\in \R$ and $t>0$ define
\begin{equation}\label{LPPtime}
L^{\lambda}_t(x)\df\sup\left\{ M_0^\lambda(z) + L((z,0),(x,t))\,:\,z\in (-\infty,x]\right\}\,,
\end{equation}
where
$$M^\lambda(z)\df\left\{\begin{array}{ll}\eta_0^\lambda((0,z])& \mbox{ for } z> 0\\
-\eta_0^\lambda((z,0]) & \mbox{ for } z\leq0\,.\end{array}\right.$$ 
Then $L_t^\lambda(x)$ is well defined (Lemma 6 from \cite{AlDi}), and it can be seen as the flux of particles through the space-time line from $(0,0)$ to $(x,t)$: 
\begin{equation}\label{LPP}
\eta^\lambda_t((x,y])=L_t^\lambda(y)-L_t^\lambda(x)\,,
\end{equation}
for all $t\geq 0$ and $x,y\in\R$ with $x<y$ (Figure 3). In the exponential LPP model we have a similar evolution given by \eqref{ExpLPP}.  

\begin{figure}
\begin{center}
\includegraphics[width=10cm, height=6cm]{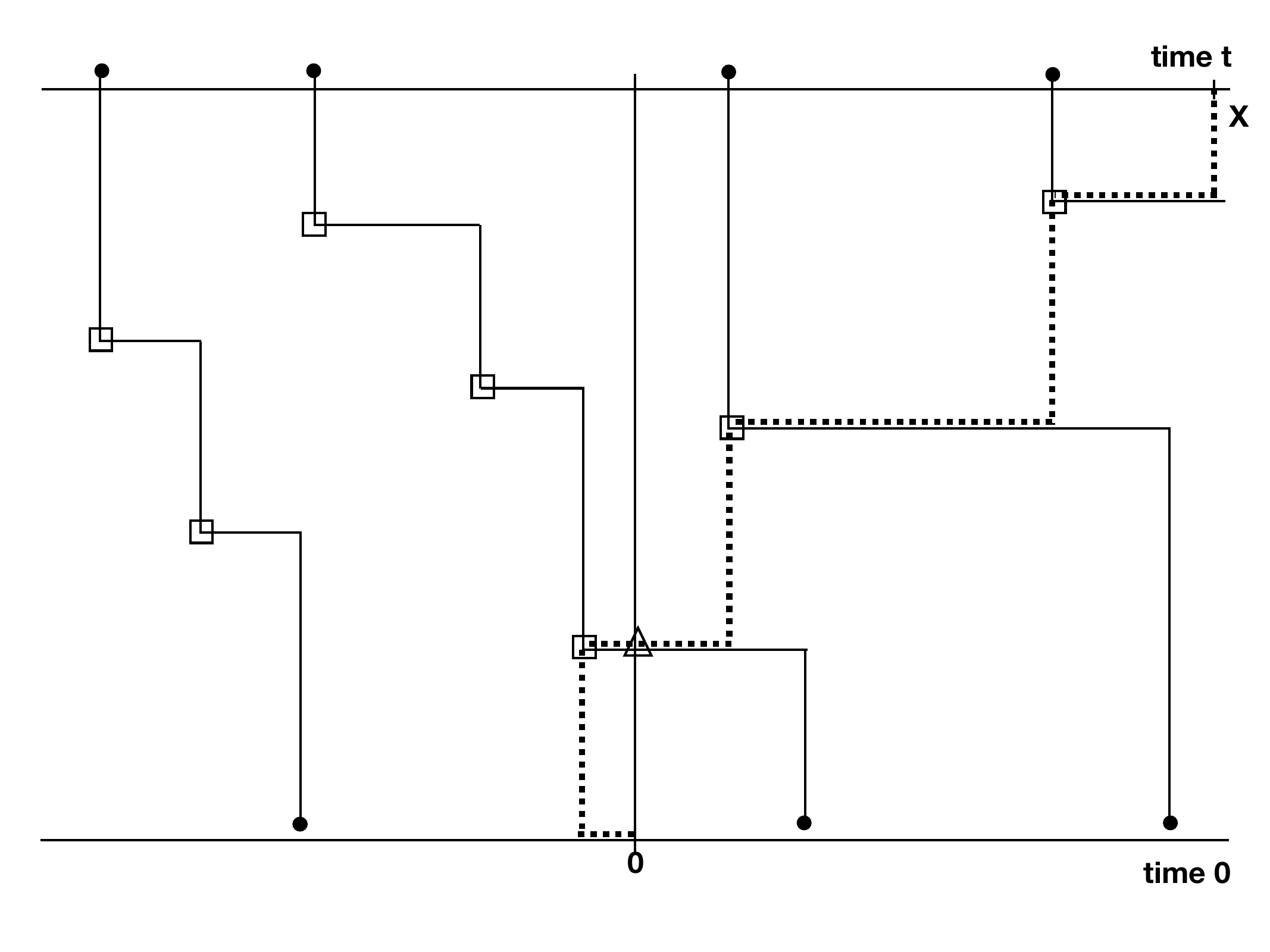}
\caption{In this picture, $L^\lambda_t(x)=3$, $L_t^\lambda(0)=1$ and $\eta^\lambda_t((0,x])=2$. A maximal path from $(0,0)$ to $(x,t)$ is drawn with dotted broken lines. In the sources ($\bullet$) and sinks ($\triangle$) formulation with $B=[0,x]\times[0,t]$, we have two sources in $[0,x]\times\{0\}$ and one sink in $\{0\}\times[0,t]$. In this setting, a  maximal path picks one sink in $\{0\}\times[0,t]$ and then follows a maximal path through the Poisson clocks ($\square$).}\label{FigLPP}
\end{center}
\end{figure}

Now we introduce an useful representation of Hammersley's process in terms of sources and sinks \cite{CaGr1}. This representation is analogue to the boundary condition with parameter $\alpha\in(0,1)$ that we have introduce in the exponential LPP model. Given a space-time box $B=[x_0,x_1]\times[t_0,t_1]$, denote by $\omega^B_v$ the marks (called sinks) generated by the flux of particles through $\{x_0\}\times(t_0,t_1)$. Then for all $(s,t]\subseteq(t_0,t_1)$ we have that  
$$\omega^B_v((s,t])\df L_t^{\lambda}(x_0)-L_s^{\lambda}(x_0)\,,$$ and $\omega^B_v$ is a Poisson point process of parameter $1/\lambda$ \cite{CaGr1}. Denote by $\omega^B_h$ the position of the particles (called sources) in $(x_0,x_1)\times\{t_0\}$. Then for all $(x,y]\subseteq(x_0,x_1]$ we have that
$$\omega^B_h((x,y])\df L_{t_0}^{\lambda}(y)-L_{t_0}^{\lambda}(x)\,, $$
and $\omega^B_h$ is a Poisson point process of parameter $\lambda$ \cite{CaGr1}. Finally let $\omega^B_b\df (x_0,x_1)\times(t_0,t_1)\cap\omega$ denote the space-time clocks within $B$ (Figure 3). Then the vector $(\omega^B_v\,,\,\omega^B_h\,,\,\omega^B_b)$, with  sinks, sources and clocks, has independent entries \cite{CaGr1}. The dynamic programming property for last-passage percolation times states that we can represent the particle trajectories inside $B$ as 
\begin{equation}\label{DynProg1}
\H^\lambda(B)\df  \left(L_t^{\lambda,B}(x)\,,\,(x,t)\in B\right)= F(\omega^B_v\,,\,\omega^B_h\,,\,\omega^B_b)\,,
\end{equation}
where $F$ is the optimization deterministic function that constructs the last-passage percolation times $L^{\lambda,B}$ inside $B$ using  $(\omega^B_v\,,\,\omega^B_h\,,\,\omega^B_b)$ \cite{CaGr1} as follows:
$$(x,t)\in B\mapsto L^{\lambda,B}_t(x)\df\sup\left\{M^{\lambda,B}(z)+L^{B}\left(z,(x,t)\right)\,:\,z\in [-(t_1-t_0),x_1-x_0]\right\}\,,$$
where 
$$L^{B}(z,(x,t))\df\begin{cases}L((z+x_0,t_0),(x,t))\,\,\,\,\,\mbox{ if }z\in(0,x_1-x_0]\\
L((x_0,t_0-z),(x,t))\,\,\mbox{ if }z\in[t_0-t_1,0]\,,
\end{cases} $$
and
$$M^{\lambda,B}(z)\df\begin{cases}M^{\lambda,B}_h(z)\df\omega^B_v\left((x_0,x_0+z]\times\{t_0\}\right)\,\mbox{ if }z\in(0,x_1-x_0]\\
M^{\lambda,B}_v(z)\df\omega^B_h\left(\{x_0\}\times(t_0,t_0-z]\right)\,\,\mbox{ if }z\in[t_0-t_1,0]\,.
\end{cases} $$
If $A$ is a box, with southwest vertex $(a,s)$, that is contained in $B$ we define 
$$\H^\lambda(B)\mid_{A}\df\left(L_t^{\lambda,B}(x)-L_s^{\lambda,B}(a)\,,\,(x,t)\in A\right)\,.$$
The dynamic programming property for last-passage percolation times also implies that 
\begin{equation}\label{DynProg2}
\mbox{ if }A\subseteq B\mbox{ then }\H^\lambda(B)\mid_{A}=\H^\lambda(A)\,.
\end{equation}

\subsection{Domination via Exit-Points}\label{DomExitPoint}
The suprema in the definition of $L^\lambda_t(x)$ is attained (possible in more than one point $z$), and one can define the right most ``exit-point'' as  
\begin{equation}\label{LPPexit}
Z^\lambda_t(x) \df\sup\left\{z\in(-\infty,x]\,:\,L^\lambda_t(x)=M^\lambda(z)+L((z,0),(x,t))\right\}\,.
\end{equation}
The picture we have in mind is that, in the last-passage percolation model with boundary $M^\lambda$, a maximal path first collects the Poisson points in the horizontal axis until the ``exit point''$Z^\lambda$, and then  collects the Poisson clocks following and increasing path. In Figure \ref{FigLPP}, the right most maximal path goes left $(Z^\lambda_t(x)<0)$ to reach the maximal path thought the Poisson clocks. It is worth to notice the following monotonicity of exit-points: 
$$Z^\lambda_t(x)\leq Z^\lambda_t(y)\mbox{ for all }y\geq x\mbox{ and }Z^\lambda_u(x)\leq Z^\lambda_t(x)\mbox{ for all }u\geq t\,.$$ 
The graphical construction allows us to run simultaneously two Hammersley's interacting particle processes with different initial profiles and using the same ``clocks'' $\omega$, known as the basic coupling. If we take $\eta^\lambda_0$ and $\eta^{\lambda'}_0$ independently, then $\eta^{\lambda'}$ is independent of $\eta^{\lambda}_0$.  The key to compare particle configurations starting from different profiles, under the basic coupling, is the following lemma. Similar versions of this lemma  can be found in \cite{CaGr2, CaPi1, CaPi2, Pi}, and the proof in our context follows the same lines. For the sake of completeness, we decide to include it in this article.
\begin{lem}\label{Comparison}
Let $a<b$ and $t\geq 0$. If $Z^{\lambda}_t(b)\leq Z^{\lambda'}_t(a)$ then
$$L^{\lambda}_t(y)-L_t^{\lambda}(x) \leq L_t^{\lambda'}(y)-L_t^{\lambda'}(x)\,,\,\mbox{ for all }x,y\in[a,b] \mbox{ with }x<y\,.$$
\end{lem}

\noindent{\bf Proof of Lemma \ref{Comparison}\,\,} 
 An increasing path $\gamma$ from $(z,s)$ to $(x,t)$ will be viewed as the lowest non-decreasing  continuous path connecting all the points, starting at $(z,s)$ and ending at $(x,t)$. In this way, we can talk about crossings with other paths or with lines. The geodesic between $(z,s)$ and $(x,t)$ is given by the lowest path (in the sense we just described) that attains the maximum in the definition of $L((z,s),(x,t))$. We will denote this geodesic by $\pi((z,s),(x,t))$. Notice that
$$L((z,s),(x,t))=L((z,s),(y,r))+L((y,r),(x,t))\,,$$
for any $(y,r)\in\pi((z,s),(x,t))$. 

Denote $z_1\df Z_t^\lambda(y)$ and $z_2\df Z_t^{\lambda'}(x)$. Let $\cc$ be a crossing between the two geodesics $\pi((z_1,0),(y,t))$ and $\pi((z_2,0),(x,t))$. Such a crossing exists  because $a<x<y<b$ and 
$$z_1=Z_t^\lambda(y)\leq Z_t^\lambda(b)\leq Z_t^{\lambda'}(a)\leq Z_t^{\lambda'}(x)=z_2\,.$$
We remark that, by superaddivity,
$$L^{\lambda'}_t(y) \geq  M^{\lambda'}(z_2) + L((z_2,0),(y,t)) \geq  M^{\lambda'}(z_2) + L((z_2,0),\cc) + L(\cc,(y,t))\,.$$
We use this, and that (since $\cc\in\pi([z_2]_0,[x]_t)$)
$$ M^{\lambda'}(z_2) + L((z_2,0),\cc)-L^{\lambda'}_t(x)= -L(\cc,(x,t))\,,$$
in the following inequality:
\begin{eqnarray*}
L^{\lambda'}_t(y) - L^{\lambda'}_t(x) & \geq & M^{\lambda'}(z_2)+L((z_2,0),\cc) + L(\cc,(y,t)) - L^{\lambda'}_t(x)\\
& = & L(\cc,(y,t)) - L(\cc,(x,t))\,.
\end{eqnarray*}
By superaddivity,
$$ - L(\cc\,,\,(x,t))\geq L^{\lambda}(\cc)-L_t^{\lambda}(x)\,,$$
and hence (since $\cc\in\pi((z_1,0),(y,t))$)
\begin{eqnarray*}
L^{\lambda'}(y,t) - L^{\lambda'}_t(x) & \geq & L(\cc,(y,t)) - L(\cc,(x,t))\\
& \geq & L(\cc,(y,t)) + L^{\lambda}(\cc)-L^{\lambda}_t(x)\\
& = & L^{\lambda}_t(y)-L^{\lambda}_t(x)\,.
\end{eqnarray*}

\hfill$\Box$\\ 

Looking at \eqref{LPP} and Lemma \ref{Comparison} we see that, to dominate $\eta^\lambda_t$ by $\eta_t^{\lambda'}$ in a finite interval $[a,b]$, one needs to control the relative positions of $Z^{\lambda}_t(b)$ and $Z^{\lambda'}_t(a)$, which can be obtained by using a large deviation bound for exit points. By translation invariance of $\eta_0^\lambda$ and $\omega$, 
\begin{equation}\label{Translation}
Z_t^\lambda(x+h)\dst Z_t(x)+h\,,
\end{equation}
and by the invariance of $\omega$ under the map $(x,t)\mapsto (\lambda x,t/\lambda)$,
\begin{equation}\label{Symmetry}
\lambda Z^\lambda_t(x)\dst  Z^1_{t/\lambda}(\lambda x)\,.
\end{equation}
The exit point $Z^\lambda_t(x)$ satisfies a law of large numbers \cite{CaGr1,CaGr2},
$$\lim_{t\to\infty}\frac{Z^\lambda_t(\lambda^2 t)}{t}\as 0\,.$$
Together with \eqref{Translation}, this says that the location of the exit-point follows, going backward in time, the characteristic direction of the system. To control large deviations, we prove the following lemma.
\begin{lem}\label{LargeDev}
 There exists an universal constant $c_0>0$ such that for all $\delta \in(0,3/4]$ 
$$\max\Big\{\P\left[Z^1_t(t)>\delta t \right]\,,\,\P\left[Z^1_t(t)<-\delta t \right]\Big\}\leq 5 e^{-c_0\delta^4 t}\mbox{ for all }t>0\,.$$
\end{lem}

\noindent{\bf Proof of Lemma \ref{LargeDev}\,\,}
By stationarity, one can think that we have Hammersley's process running for $t\in\R$. The flux of particles through the vertical axis $(L^\lambda_t(0)\,,\,t\in\R)$, seen as a point process, is also Poisson with parameter $1/\lambda$. By using reflection in the diagonal, 
$$\P\left[Z^\lambda_t(x)<0\right]\leq \P\left[Z^{1/\lambda}_x(t)>0\right] \,.$$  
By combining this with translation invariance \eqref{Translation}, we get that 
$$\P\left[Z^1_t(t)<-u \right]=\P\left[Z^1_t(t+u)<0 \right]\\ 
\leq \P\left[Z^1_{t+u}(t)>0 \right]=\P\left[Z^1_{t+u}(t+u)>u \right]\,,$$
and to obtain the desired bound we only need to deal with $\P\left[Z^1_t(t)>\delta t \right]$. For that, we parallel the arguments developed in \cite{CaGr2} (see Theorem 4.2 there) in which they prove $2/3$ asymptotics for $Z^\lambda$, which also works for large deviation bounds.  

Let $\lambda>1$ and take $M^\lambda$ a Poisson process of parameter $\lambda$ obtained by thickening of $M^1$. By definition of $L^\lambda$,
$$L\left((z,0),(t,t)\right)\leq L_t^\lambda(t)-M^\lambda(z)\,,\mbox{ for all }z\leq t\,,$$
and hence 
\begin{eqnarray*}
\P\left[ Z^1_t(t)>u\right]&=&\P\left[ \exists \,z\in( u,t]\,:\,M^1(z)+L\left((z,0),(t,t)\right)=L^1_t(t)\right]\\
&\leq &\P\left[ \exists \,z\in( u,t]\,:\,M^1(z)+L_t^\lambda(t)-M^\lambda(z)\geq L^1_t(t)\right]\\
&=&\P\left[ \exists \,z\in( u,t]\,:\,M^{\lambda}(z)-M^1(z)\leq L^\lambda_t(t)-L^1_t(t)\right]
\end{eqnarray*}
By construction, $M^{\lambda}(z)-M^1(z)$ is increasing as a function of $z$, and $M^{\lambda}-M^1\dst M^{\lambda-1}$. Thus
$$\P\left[ Z^1_t(t)>u\right]\leq \P\left[  M^{\lambda}(u)-M^1(u)\leq L^\lambda_t(t)-L^1_t(t)\right]\,, $$ 
and to obtain a useful bound we maximize 
$$\E\left[M^{\lambda-1}(u) \right]-\E\left[ L^\lambda_t(t)-L^1_t(t) \right]=(\lambda-1)u-t\left(\lambda+\frac{1}{\lambda}-2\right)\,,$$
as a function of $\lambda>1$, and choose
$$\lambda_{u/t}=\left(1-u/t\right)^{-1/2}>1\,.$$
For all $u\leq \frac{3}{4} t$ we have that 
$$\E\left[M^{\lambda_{u/t}-1}(u) \right]-\E\left[ L^{\lambda_{u/t}}_t(t)-L^1_t(t) \right]\geq \frac{u^2}{4t}\,,$$
and hence
\begin{eqnarray*}
\P\left[ Z^1_t(t)>u\right]&\leq& \P\left[  M^{\lambda_{u/t}}(u)-M^{1}(u)\leq L^{\lambda_{u/t}}_t(t)-L^1_t(t)\right]\\
&\leq &\P\left[  M^{\lambda_{u/t}-1}(u)\leq \E\left[M^{\lambda_{u/t}-1}(u)\right]-u^2/(8t)\right]\\
&+&\P\left[  L^{\lambda_{u/t}}_t(t)-L^1_t(t) \geq \E\left[M^{\lambda_{u/t}-1}(u)\right]-u^2/(8t)\right]\\
&\leq & \P\left[  M^{\lambda_{u/t}-1}(u)\leq \E\left[M^{\lambda_{u/t}-1}(u)\right]-u^2/(8t)\right]\\
&+&\P\left[  L^{\lambda_{u/t}}_t(t)-L^1_t(t) \geq \E\left[L^{\lambda_{u/t}}_t(t)-L^1_t(t) \right]+u^2/(8t)\right]\,.
\end{eqnarray*}
Recall that $u=\delta t$ for $\delta\in(0,3/4]$, and   
$$\lambda_{u/t}=\lambda_\delta=(1-\delta)^{-1/2}\in(1,2]\,\,\mbox{ and }\,\,0<\frac{\lambda_\delta-1}{\delta}\leq 2\,.$$
In what follows, we use that if $\Poi^{\lambda}$ denotes a Poisson random variable with parameter $\lambda>0$ then 
$$\max\big\{\P\left[\Poi^{\lambda}\geq \lambda +x\right],\P\left[\Poi^{\lambda}\leq \lambda -x\right]\big\}\leq e^{-\frac{x^2}{2(x+\lambda)}}\,,\,\forall \,x>0\,.$$ 
See Lemma \ref{PoissonConcentration} in the Appendix \ref{Poisson}. Hence 
\begin{multline*}
\P\left[  M^{\lambda_{\delta}-1}(u)\leq \E\left[M^{\lambda_{\delta}-1}(u)\right]-u^2/(8t)\right] \leq \exp\left( -\frac{1}{2}\left(\frac{\frac{\delta^4}{8^2}t^2}{\frac{\delta^2}{8} t+(\lambda_\delta-1)\delta t}\right)\right)\\=\exp\left(-\left(\frac{1}{1+\frac{ 8(\lambda_\delta-1)}{\delta}} \right) \frac{\delta^2}{16} t\right)
\leq \exp\left(-\left(\frac{1}{1+16} \right) \frac{\delta^2}{16} t\right) =\exp\left(-\frac{\delta^2}{272} t\right)\,.
\end{multline*} 
For the other term, notice that  
\begin{multline*}
P\left[ L^{\lambda_{\delta}}_t(t)-L^1_t(t) \geq \E\left[L^{\lambda_{\delta}}_t(t)-L^1_t(t) \right]+u^2/(8t)\right]\\
= \P\left[ \left( L^{\lambda_\delta}_t(t)- \E\left[L^{\lambda_\delta}_t(t)\right]\right)+ \left(\E\left[L^1_t(t) \right] -L^1_t(t)\right) \geq \frac{\delta^2}{8} t\right]\\
\leq\P\left[ L^{\lambda_\delta}_t(t)- \E\left[L^{\lambda_\delta}_t(t)\right] \geq \frac{\delta^2}{16} t\right]+\P\left[L^1_t(t)-\E\left[L^1_t(t) \right]  \leq -\frac{\delta^2}{16} t\right]\,.
\end{multline*}
For arbitrary $\lambda>0$ and $x,t>0$, 
$$L^{\lambda}_t(x)=L^{\lambda}_t(0)+\left(L^{\lambda}_t(x)-L^{\lambda}_t(0)\right)\,,$$
where $L^{\lambda}_t(0)\dst\Poi^{t/\lambda}$ and $L^{\lambda}_t(x)-L^{\lambda}_t(0)\dst\Poi^{\lambda x}$. Thus
$$\P\left[ L^{\lambda_\delta}_t(t)- \E\left[L^{\lambda_\delta}_t(t)\right] \geq \frac{\delta^2}{16} t\right]\leq 
\P\left[ \Poi^{t/\lambda_\delta} \geq\frac{t}{\lambda_\delta}+\frac{\delta^2}{32} t\right]+\P\left[\Poi^{\lambda_\delta t}\geq \lambda_\delta t+\frac{\delta^2}{32} t\right]\,,$$
and 
$$\P\left[ L^{1}_t(t)- \E\left[L^{1}_t(t)\right] \leq-\frac{\delta^2}{16} t\right]\leq 
2\P\left[ \Poi^{t} \leq t-\frac{\delta^2}{32} t\right]\,.$$
Since $\frac{\delta^2}{32}+1\leq\frac{9}{512}+1=\frac{521}{512}$,
$$\P\left[ \Poi^{t} \leq t-\frac{\delta^2}{32} t\right]\leq \exp\left(-\frac{1}{2}\left(\frac{\frac{\delta^4}{32^2}}{\frac{\delta^2}{32}+1}\right)t \right) \leq  \exp\left(-\frac{1}{2}\left(\frac{\frac{1}{1024}}{\frac{521}{512}}\right)\delta^4t \right)\,.$$
An analog bound can be obtained for $\lambda_\delta\in(1,2]$, which concludes the proof.

\hfill$\Box$\\

By combining Lemma \ref{LargeDev} together with properties \eqref{Translation} and \eqref{Symmetry}, we are able to obtain a large deviation bound for the location of the exit point for any $x \in \mathbb{R}$ and $t>0$, at any density $\lambda>0$:
\begin{equation}\label{ld}
 \max\Big\{\P\left[Z^\lambda_t(x)> x - \frac{t}{\lambda^2} + \delta \frac{t}{\lambda^2} \right]\,,\,\P\left[Z^\lambda_t(x)<x - \frac{t}{\lambda^2} - \delta \frac{t}{\lambda^2} \right]\Big\}\leq 5 e^{-c_0\delta^4 t/\lambda}\,.   
\end{equation}
Similarly to \cite{BaTe,HiHoSiSoTe}, one of the main ingredient to prove Theorem \ref{Sprinkling} is the construction of a coupling that produces domination between particle configurations with high probability, which is achieved by using Lemma \ref{Comparison} and \eqref{ld}.   
\begin{lem}\label{Coupling}
Let $\lambda<\lambda'$ and assume that $\epsilon\df\left(\lambda^{-2}-\lambda'^{-2}\right)\in\left(0,\frac{3}{2}\lambda'^{-2}\right)$.  Under the basic coupling $(\eta^{\lambda},\eta^{\lambda'})$, with $\eta_0^{\lambda}$ independent of  $\eta_0^{\lambda'}$, we have that $\eta^{\lambda'}$ is independent of $\eta^{\lambda}_0$. Furthermore, for all $a\geq b$ and $s\geq 0$ we have that 
$$\P\left[\forall\,u\in[t,t+s]\,,\,\forall\,(x,y]\subseteq(a,b]\,,\,\eta^{\lambda}_u((x,y])\leq\eta^{\lambda'}_u((x,y]) \right]\geq 1- 10 e^{-\frac{c_0\lambda^7}{4^4}\epsilon^4 t}\,,$$
for all $t>4\epsilon^{-1}\left(b-a+s\lambda^{-2}\right)$.
\end{lem}

\noindent{\bf Proof of Lemma \ref{Coupling}\,\,} 
Notice that $L^{\lambda}$ and $L^{\lambda'}$ are naturally coupled since they are functions of the respective initial conditions and of the same two-dimensional Poisson point process $\omega$ (basic coupling). By monotonicity of $Z^\lambda$, \eqref{LPP} and Lemma \ref{Comparison} we have that
$$ Z^{\lambda}_t(b)\leq Z^{\lambda'}_{t+s}(a)\,\implies\,\forall\,u\in[t,t+s]\,\forall\,(x,y]\subseteq(a,b]\,,\,\eta^{\lambda}_t((x,y])\leq\eta^{\lambda'}_t((x,y])\,,$$
and it suffices to prove that 
$$\P\left[Z^{\lambda}_t(b)> Z^{\lambda'}_{t+s}(a) \right]\leq 10 e^{-\frac{c_0\lambda^7}{4^4}\epsilon^4 t}\,.$$
Since for all $m\in\R$
$$\left\{Z^{\lambda}_t(b)\leq m \right\}\cap\left\{m\leq Z^{\lambda'}_{t+s}(a) \right\} \subseteq \left\{Z^{\lambda}_t(b)\leq Z^{\lambda'}_{t+s}(a) \right\}\,,$$
we have that 
$$\P\left[Z^{\lambda}_t(b)> Z^{\lambda'}_{t+s}(a) \right]\leq \P\left[ Z^{\lambda'}_{t+s}(a)< m\right]+\P\left[Z^{\lambda}_t(b)> m \right] \,.$$

Define $m\df a-\frac{\lambda^{-2}+\lambda'^{-2}}{2}(t+s)$ and note that 
$$m=a-\lambda'^{-2}(t+s)- \lambda'^{-2}\delta'(t+s)\,,\mbox{ where }\delta'=\lambda'^{2}\frac{\epsilon}{2}<\frac{3}{4}\,.$$
By \eqref{ld},  
$$\P\left[ Z^{\lambda'}_{t+s}(a)< m\right]\leq 5 e^{-\frac{c_0}{\lambda'}\delta'^4 (t+s)}\leq 5 e^{-\frac{c_0\lambda'^7}{2^4}\epsilon^4 t}\leq 5 e^{-\frac{c_0\lambda^7}{4^4}\epsilon^4 t}\,. $$
To bound the other term, we take $\delta=\lambda^2\frac{\epsilon}{4}<\frac{3}{8}$ and note that 
$$\left(\delta\frac{t}{\lambda^2}+(b-a)+\frac{s}{\lambda^2}-(t+s)\frac{\epsilon}{2}\right)<0\,\iff\,t\left(\frac{\epsilon}{2}-\frac{\delta}{\lambda^2}\right)>s\left(\frac{\lambda^{-2}+\lambda'^{-2}}{2}\right)+b-a\,,$$
which is true for $t>4\epsilon^{-1}\left(\frac{s}{\lambda^2}+b-a\right)$. Thus, we can write 
$$m=b-\frac{t}{\lambda^2}+\delta\frac{t}{\lambda^2}-\left(\delta\frac{t}{\lambda^2}+(b-a)+\frac{s}{\lambda^2}-(t+s)\frac{\epsilon}{2}\right) > b-\frac{t}{\lambda^2}+\delta\frac{t}{\lambda^2}\,.$$
By \eqref{ld}, this implies that 
$$\P\left[ Z^{\lambda}_t(b)> m\right]\leq \P\left[ Z^{\lambda}_t(b)>b-\frac{t}{\lambda^2}+\delta\frac{t}{\lambda^2} \right]\leq 5e^{-c_0\frac{\delta^4}{\lambda}t}\leq 5 e^{-\frac{c_0\lambda^7}{4^4}\epsilon^4 t} \,,$$
for $t>4\epsilon^{-1}\left(\frac{s}{\lambda^2}+b-a\right)$.

\hfill$\Box$\\ 

\subsection{Lateral Decoupling}
Let $d_h,b_1,b_2,s_1,s_2>0$ and $\delta_v\geq -s_1$. Denote  
\begin{equation}\label{LeftBelow}
B_1=[-b_1,0]\times[0,s_1]\,\mbox{ and }\,B_2=[d_h,d_h+b_2]\times[\delta_v+s_1,\delta_v+s_1+s_2]\,.
\end{equation}
We note that $d_h$ gives the horizontal distance between $B_1$ and $B_2$, while for the vertical distance we have that $d_v=\max\{\delta_v,0\}$. Now consider the space-time box 
$$B_{1,2}\df [0,d_h+b_2]\times [0,\delta_v+s_1+s_2]\,.$$
Recall \eqref{DynProg1} and construct last-passage percolation times with no sinks\footnote{We represent the identically zero measure by $\emptyset$.} by setting
$$\widetilde{\H}^\lambda(B_{1,2})\df F(\emptyset\,,\,\omega^{B_{1,2}}_h\,,\,\omega^{B_{1,2}}_b)\,\mbox{ for }(x,t)\in B_{1,2}\,.$$
Since $(\emptyset\,,\,\omega^{B_{1,2}}_h\,,\,\omega^{B_{1,2}}_b)$ is independent of $(\omega^{B_1}_v\,,\,\omega^{B_{1}}_h\,,\,\omega^{B_{1}}_b)$, by \eqref{DynProg2}, we have that 
\begin{equation}\label{Indep1}
\H^\lambda(B_1)\mbox{ is independent of }\widetilde{\H}^\lambda(B_{1,2})\mid_{B_2}=\widetilde{\H}^\lambda(B_2)\,.
\end{equation}
Furthermore,
\begin{equation}\label{Tilde1}
Z^\lambda_{\delta_v+s_1+s_2}(d_h)\geq 0\,\implies\,\H^\lambda(B_2)=\widetilde{\H}^\lambda(B_2)\,.
\end{equation}
To justify \eqref{Tilde1}, we begin by noticing that $\widetilde{L}^\lambda(x,t)\leq L^\lambda(x,t)$ for all $(x,t)\in B_{1,2}$ since an admissible path for the environment $(\emptyset\,,\,\omega^{B_{1,2}}_h\,,\,\omega^{B_{1,2}}_b)$ is also an admissible path for the environment $(\omega^{B_{1,2}}_v\,,\,\omega^{B_{1,2}}_h\,,\,\omega^{B_{1,2}}_b)$. By monotonicity, if $Z_{\delta_v+s_1+s_2}^\lambda(d_h)\geq 0$ then $Z_t^\lambda(x)\geq 0$ for all $(x,t)\in B_2$, and a maximal path for $(\omega^{B_{1,2}}_v\,,\,\omega^{B_{1,2}}_h\,,\,\omega^{B_{1,2}}_b)$, with $Z_t^\lambda(x)\geq 0$, will also be a maximal path for $(\emptyset\,,\,\omega^{B_{1,2}}_h\,,\,\omega^{B_{1,2}}_b)$, which implies that $\widetilde{L}^\lambda(x,t)= L^\lambda(x,t)$. Together with \eqref{DynProg2}, this implies that 
$$\widetilde{\H}^\lambda(B_2)=\widetilde{\H}^\lambda(B_{1,2})\mid_{B_2}=\H^\lambda(B_{1,2})\mid_{B_2}=H^\lambda(B_2)\,.$$

Now if we denote
\begin{equation}\label{RightBelow}
B_1=[d_h,d_h+b_1]\times[0,s_1]\,\mbox{ and }\,B_2=[-b_2,0]\times [\delta_v+s_1,\delta_v+s_1+s_2]\,,
\end{equation}
and let
$$B_{1,2}\df [0,d_h+b_1]\times [0,\delta_v+s_1+s_2]\,,$$
then $\widetilde\H^\lambda(B_{1,2})$ is a function of the sources $\omega_h^{B_{1,2}}$ and the clocks $\omega_b^{B_{1,2}}$, while $\H^\lambda(B_2)$ can be seen as a function of the sources to the left of the origin and the clocks on the second quadrant $\{(x,t)\,:\,x<0\,,\,t> 0\}$. Therefore,
\begin{equation}\label{Indep2}
\H^\lambda(B_2)\mbox{ is independent of }\widetilde{\H}^\lambda(B_{1,2})\mid_{B_1}=\widetilde{\H}^\lambda(B_1)\,.
\end{equation}
In addition, 
\begin{equation}\label{Tilde2}
Z_{s_1}(d_h)\geq 0\,\implies\,\widetilde\H^\lambda(B_{1})=\widetilde\H^\lambda(B_{1,2})\mid_{B_1}=\H^\lambda(B_1)\,.
\end{equation}
The proof of \eqref{Tilde2} follows the same lines of the proof of \eqref{Tilde1}.

Notice that if $d_h\geq c_1\left(d_v+\per(B_1)+\per(B_2)\right)$ then $d_h\geq c_1 \max\{\delta_v+s_1+s_2\,,\, s_1\}$, and by monotonicity of the exit-point, this implies that 
\begin{equation}\label{LateralIneq}
\max\left\{\P\left[ Z^\lambda_{\delta_v+s_1+s_2}(d_h)<0 \right]\,,\,\P\left[Z_{s_1}(d_h)<0\right]\right\}\leq\P\left[Z_{d_h/c_1}(d_h)<0\right]\,.
\end{equation}
Thus, to prove asymptotic independence between boxes that are horizontally far apart (lateral decoupling) we use the next lemma. 
\begin{lem}\label{Lateral}
There exist constants $c_i=c_i(\lambda)>0$ for $i=1,2$ such that  
$$ \P\left[ Z^\lambda_{d_h/c_1}(d_h)<0 \right]\leq 5 e^{-c_2 d_h}\,,\mbox{ for all }d_h>0\,.$$
\end{lem}

\noindent{\bf Proof of Lemma \ref{Lateral}\,\,} Apply \eqref{ld} with  $x=d_h$, $t = d_h/c_1$ and $\delta = 3/4$,
$$\P\left[Z_{d_h/c_1}(d_h)< d_h - \frac{d_h}{c_1 \lambda^2} - \frac{3 d_h}{4 c_1 \lambda^2}\right] = \P\left[Z_{d_h/c_1}(d_h)< d_h \left( 1 - \frac{7}{4 c_1 \lambda^2} \right) \right] \leq 5 e^{-c_0\left(\frac{3}{4}\right)^4 \frac{d_h}{c_1\lambda}}\,,$$
and take $c_1$  such that  $1 - \frac{7}{4 c_1 \lambda^2} >0$. This finishes the proof of the lemma.

\hfill$\Box$\\ 

\subsection{Proof of Theorem \ref{Sprinkling}}

Without loss of generality, we assume that $B_1$ is below $B_2$. We proceed with the proof by first supposing that  
$$d_h\geq c_1\left(d_v+\per(B_1)+\per(B_2)\right)\,.$$
If $B_1$ is to the left of $B_2$, by translation invariance we can assume \eqref{LeftBelow}. Given a function $f_2$ of $\H^\lambda(B_2)$, we denote $\widetilde{f}_2$ the same function but now applied to  $\widetilde{\H}^\lambda(B_2)$. Since $f_2,\widetilde f_2\in[0,1]$, by \eqref{Tilde1} and \eqref{LateralIneq},
$$\E^\lambda\left[|f_2-\widetilde{f}_2|\right]= \E^\lambda\left[|f_2-\widetilde{f}_2|\1_{\{Z^\lambda_{\delta_v+s_1+s_2}(d_h)<0\}}\right]\leq \P\left[Z^\lambda_{d_h/c_1}(d_h)<0\right] \,.$$
By \eqref{Indep1}, if $f_1$ is a function of $\H^\lambda(B_1)$ then $f_1$ is independent of $\widetilde f_2$ and 
$$\E^\lambda\left[f_1f_2\right]=\E^\lambda\left[f_1\right]\E^\lambda\left[f_2\right]+\E^\lambda\left[f_1\right]\left(\E^\lambda\left[\widetilde{f}_2\right]-\E^\lambda\left[f_2\right]\right)+\E^\lambda\left[f_1(f_2-\widetilde{f}_2)\right]\,,$$
which implies that 
$$\Big|\E^\lambda\left[f_1f_2\right]-\E^\lambda\left[f_1\right]\E^\lambda\left[f_2\right]\Big|\leq 2 \E^\lambda\left[|f_2-\widetilde{f}_2|\right]\leq 2\P\left[Z^\lambda_{d_h/c_1}(d_h)<0\right]\,.$$
By Lemma \ref{Lateral}, 
$$\Big|\E^\lambda\left[f_1f_2\right]-\E^\lambda\left[f_1\right]\E^\lambda\left[f_2\right]\Big|\leq 4 e^{-c_2 d_h}\leq 4 e^{-\frac{c_2}{1+c_1^{-1}}d}\,,$$
since $d=d_h+d_v\implies d\leq d_h+ d_h/c_1-\per(B_1)-\per(B_2) \leq (1+c_1^{-1})d_h$. If $B_1$ is to the right $B_2$, by translation invariance we can assume \eqref{RightBelow}. Given a function $f_1$ of $\H^\lambda(B_1)$, we denote $\widetilde{f}_1$ the same function but now applied to  $\widetilde{\H}^\lambda(B_1)$. By using \eqref{Indep2}, \eqref{Tilde2} and \eqref{LateralIneq}, one can go on as before, show that 
$$\Big|\E^\lambda\left[f_1f_2\right]-\E^\lambda\left[f_1\right]\E^\lambda\left[f_2\right]\Big|\leq 2 \E^\lambda\left[|f_1-\widetilde{f}_1|\right]\leq 2\P\left[Z^\lambda_{d_h/c_1}(d_h)<0\right]\,,$$
and apply Lemma \ref{Lateral} to conclude the same upper bound. By looking at the definition of $c_1$ and $c_2$ in Lemma \ref{Lateral}, one can see that $c_1$ and $c_2$ can be chosen so that the upper bound holds for all $\lambda\in[\lambda_1 ,\lambda_2]\subseteq (0,\infty)$.     

Now, for $d_h < c_1 (d_v+\per(B_1)+\per(B_2))$ and $d=d_v+d_h$ we have that 
$$d< (1+c_1)d_v+c_1(\per(B_1)+\per(B_2))\,\iff\,d_v>\frac{1}{1+c_1}\left(d-c_1(\per(B_1)+\per(B_2))\right)\,.$$
Let $\epsilon<\min\left\{\lambda'^{-2}\frac{3}{2},\frac{3}{2}\right\}$, take 
$$C_\lambda= 2(1+c_1)\left(\max\left\{1,\lambda^{-2}\right\}+\frac{3c_1}{4(1+c_1)}\right)\,,$$
and assume that $d\geq \epsilon^{-1}C_\lambda\left(\per(B_1)+\per(B_2)\right)$. Therefore,
\begin{multline*}
d>2\epsilon^{-1}(1+c_1)\left(\max\left\{1,\lambda^{-2}\right\}+\frac{c_1\epsilon}{2(1+c_1)}\right) \left(\per(B_1)+\per(B_2)\right) \\
=2\epsilon^{-1}(1+c_1)\left(\max\left\{1,\lambda^{-2}\right\}\right)\left(\per(B_1)+\per(B_2)\right)+c_1\left(\per(B_1)+\per(B_2)\right) \,,
\end{multline*}
and hence,
$$d_v>\frac{1}{1+c_1}\left(d-c_1(\per(B_1)+\per(B_2))\right)\geq 2\epsilon^{-1}\left(\max\left\{1,\lambda^{-2}\right\}\right)\left(\per(B_1)+\per(B_2)\right)\,.$$
By translation invariance again, we can assume that $B_2=[a,b]\times [t,t+s]$  and $t=d_v$, so that the upper side of $B_1$ intersects the horizontal axis $\R\times\{0\}$. Then 
$$\per(B_2)=2(b-a)+2s\,\implies\,t=d_v\geq 4\epsilon^{-1} \max\{1,\lambda^{-2}\}\left((b-a)+s\right)\geq 4\epsilon^{-1}\left(b-a+s\lambda^{-2}\right) \,.$$ 
By Lemma \ref{Coupling}, if we denote 
$$E=E_{t,s,a,b}\df\left\{\forall\,u\in[t,t+s]\,,\,\forall\,(x,y]\subseteq(a,b]\,,\,\eta^{\lambda}_u((x,y])\leq\eta^{\lambda'}_u((x,y]) \right\}\,,$$
then 
$$\P\left[E^c\right]\leq 10 e^{-\frac{c_0\lambda^7}{4^4}\epsilon^4 t}\,.$$
Let $\cF \df\sigma(\eta^{\lambda}_u\,,\,u\leq 0)$. By Markov's property,
$$\E^\lambda\left[f_1 f_2\right]=\E^\lambda\left[\E\left[f_1 f_2\mid\cF\right]\right]=\E^\lambda\left[f_1\E\left[ f_2\mid\cF\right]\right]=\E^\lambda\left[f_1\E\left[ f_2\mid\eta^{\lambda}_0\right]\right]\,.$$
Recall that $f_i\in [0,1]$, that $f_i$ is non-decreasing and that, according to Lemma \ref{Coupling}, $\eta^{\lambda'}$ is independent of $\eta^\lambda_0$. If we denote $\E^\lambda\left[f\right]=\E\left[f(\eta^\lambda)\right]$ then
\begin{eqnarray*}
\E\left[f_1(\eta^\lambda)\E\left[ f_2(\eta^\lambda)\mid\eta^{\lambda}_0\right]\right]&=&\E\left[f_1(\eta^\lambda)\E\left[ f_2(\eta^\lambda)\1_{E}\mid\eta^{\lambda}_0\right]\right]+\E\left[f_1(\eta^\lambda)\E\left[ f_2(\eta^\lambda)\1_{E^c}\mid\eta^{\lambda}_0\right]\right]\\
&\leq &\E\left[f_1(\eta^\lambda)\E\left[ f_2(\eta^{\lambda'})\1_{E}\mid\eta^{\lambda}_0\right]\right]+\E\left[f_1(\eta^\lambda)\E\left[ f_2(\eta^\lambda)\1_{E^c}\mid\eta^{\lambda}_0\right]\right]\\
&\leq & \E^\lambda\left[f_1\right]\E^{\lambda'}\left[ f_2\right]+ 10 e^{-\frac{c_0\lambda^7}{4^4}\epsilon^4 t}\,.
\end{eqnarray*}
To conclude, if $d\geq C_\lambda \epsilon^{-1}(\per(B_1)+\per(B_2))$ and $d_h < c_1 (d_v+\per(B_1)+\per(B_2))$, we have that
$$
d=d_h+d_v<(1+c_1)d_v+c_1\left(\per(B_1)+\per(B_2)\right)\\
\leq (1+c_1)d_v+\frac{c_1}{C_\lambda} \epsilon d\,,$$
which implies that $(1-\frac{c_1}{C_\lambda}\epsilon)d<(1+c_1)d_v$. If we take $\epsilon<\min\left\{\lambda'^{-2}\frac{3}{2},\frac{3}{2},\frac{C_\lambda}{2c_1}\right\}$ then $d<2(1+c_1)d_v=2(1+c_1)t$ and 
$$e^{-\frac{c_0\lambda^7}{4^4}\epsilon^4 t}\leq e^{-\frac{c_0\lambda^7}{4^{9/2}(1+c_1)}\epsilon^4 d}\,,$$
which finishes the proof of Theorem \ref{Sprinkling}. It should be clear that one can choose the constants so that the condition for the inequality and the inequality itself hold for all $\lambda,\lambda'\in[\lambda_1,\lambda_2]\subseteq (0,\infty)$.

\subsection{Proof of Theorem \ref{SprinklingLPP}}
The proof is similar and follows the same geometric approach based on the control of exit points. By continuity of the exponential environment there exists a unique maximal path (or geodesic) $\pi_{kn}$ such that $G(k,n)=\sum _{(i,j)\in \pi_{kn}}w_{ij}$. As in \eqref{LPPtime} and \eqref{LPPexit}, one can define the exit point $X^\alpha_n(k)$ to represent the last point of this path on either the i-axis or the j-axis or, equivalently, the a.s. unique point such that 
$$G(k,n)=U(X^\alpha_n(k))+A(X^\alpha_n)\,,$$
where $U$ represent the weights collected on either the i-axis or the j-axis, and $A$ represents the last-passage percolation time restricted to the $\Exp(1)$ weights in the bulk. For the exponential LPP model, by continuity of the weights, we do not need to take the right most argument that maximizes this sum. Once we have that, the analogue to Lemma \ref{Comparison} follows by the same argument. Hence, in order to prove the exponential LPP version of Lemma \ref{Coupling}, we only need to have good control of the exit point location. For the exponential model we do not have invariance under linear maps that preserve volume, nevertheless given $I=[\alpha_1,\alpha_2]\subseteq(0,1)$ one can find a constant $c_I>0$ such that for all $\alpha\in I$ it holds that  
$$\max\Big\{\P\left[X^{\alpha}_n(\lfloor\beta n\rfloor)>\epsilon n \right]\,,\,\P\left[X^{\alpha}_n(\lfloor\beta n\rfloor)<-\epsilon n \right]\Big\}\leq 5 e^{-c_I\epsilon^4 n}\mbox{ for all }n\geq 1\,.$$
Recall that $\beta=\beta(\alpha)=(1-\alpha)^{-2}\alpha^2$ is the characteristic speed of the exponential LPP model with parameter $\alpha\in(0,1)$. On one hand, we do not believe that this large deviation bound is sharp, and one can possibly improve it by using Theorem 3.5 of \cite{SeSh} (see also \cite{EmJaSe}) for the tail of the exit point distribution under the $n^{2/3}$ scaling. On the other hand, it is enough for our purposes and can be obtained mutatis mutandis by paralleling the same lines as the ones for the Poisson model (Lemma \ref{LargeDev}). In order to produce independence by modifying the boundary condition (using Lemma 4.2 of \cite{BaCaSe}), as in \eqref{Indep1} and \eqref{Indep2}, one introduce null weights along the vertical axis, and notice that the analogue of \eqref{Tilde1} \eqref{Tilde2} still holds in  the lattice LPP context. To conclude lateral decoupling, one can adapt the arguments in the proof of Lemma \ref{Lateral} to the exponential model, or use sharper results such as Corollary 3.6 of \cite{SeSh}. Once all these ingredients are gathered, the proof of Theorem 4 follows the same recipe used in the proof of Theorem 1.      
\section{Renormalization Scheme for the Detection Problem}\label{Renorm}
 
The proof of Theorem \ref{Detection} uses a multiscale renormalization scheme developed by Baldasso and Texeira in \cite{BaTe}, that requires a suitable decay of correlation. If the perimeter of a box is of order $l$ and $\epsilon\sim l^{-\delta}$, then $\epsilon^4l \sim l^{1-4\delta}$ and $\epsilon^{-1}  l\sim l^{1+\delta}$. In order to use the renormalization scheme developed in \cite{BaTe}, it is enough to have that $1-4\delta>0$, to get a summable upper bound (Theorem \ref{Sprinkling}), and $1+\delta<3/2$, to ensure that boxes are far apart (see (2.24) and Lemma 2.7 \cite{BaTe}). With this in mind, $\delta<1/4$ will suffice for our purposes. We start by setting the details of the collection of oriented paths that will be used.

\begin{Def}\label{DefPaths}
Fix $R\geq 3$ and let $\cS_R$ be the collection of functions $f:\bN\to \ZZ^2$ that satisfies 
$$f(n+1)-f(n) = (r,1)\,\mbox{ for some }|r|\leq R\,.$$
\end{Def}
For $f\in\cS_R$ define $\tilde{f}:\R_+\to\R^2$ as the linear interpolation of $f$:
$$\tilde{f}(t)\df (1+\lfloor t\rfloor -t)f(\lfloor t\rfloor)+(t-\lfloor t\rfloor)f(\lfloor t\rfloor+1)\,.$$
We say that $f,g\in\cS_R$ can be concatenated when there exists $s,t\in\R_+$ such that $\tilde{f}(t)=\tilde{g}(s)=z_0$ for some $z_0 \in \ZZ^2$, and define a concatenation $h:\bN\to \ZZ^2$ of $f$ and $g$ at $z_0$ as    
$$h(n)=\begin{cases}f(n)\,\,\,\,\,\,\,\,\,\,\,\,\,\,\,\,\,\,\,\,\,\,\,\,\,\,\,\,\,\,\,\,\mbox{ if }n\geq s\,,\\ g(\lfloor t\rfloor-\lfloor s\rfloor+ n)\,\,\mbox{ if }n>s\,.\end{cases}$$
We refer the elements of $\cS_R$ as (oriented) \textit{paths}. Notice that the collection of paths $\cS_R$ is closed under concatenation in the sense that $h\in\cS_R$ for every $f$ and $g \in \cS_R$ that can be concatenated.
 
For the multiscale renormalization scheme, we follow \cite{BaTe} and define a sequence of scales 
$$l_0=10^{100}\,,\,l_{k+1}=\lfloor l_k^ {1/2}\rfloor l_k\,\mbox{ and }\,L_k=\left\lfloor\left(\frac{3}{2}+\frac{1}{k}\right)l_k\right\rfloor\,.$$
Then 
$$\frac{l_k^{3/2}}{2}\leq l_{k+1}\leq l^{3/2}_k\,\mbox{ and }\,l_k\leq L_k\leq 2l_k\,,$$
if $k$ is large enough. We also have a sequence of sets 
$$A_k\df [0,l_k]\times[0,L_k]\cup [l_k,l_k+L_k]\times[L_k,l_k+L_k]\,,$$ 
and a sequence of associated boxes $B_k\df [0,l_k+L_k]\times[0,l_k+L_k]$. We say that a path $f\in \cS_R$ is a crossing of $A_k$ if there exists $T_f>0$ such that 
\begin{itemize}
\item $f(0) \in[0,l_k] \times \{0\}$;
\item $\tilde f(T_f)\in \{l_k+L_k\}\times [L_k,l_k+L_k]$;
\item $\tilde{f}(t)\in A_k$ for all $t\in[0,T_f]$.
\end{itemize}

Given a probability space $(\Omega,\P,\cF)$, a site percolation model is based on a random element $X:\Omega\to\{0,1\}^{\ZZ^2}$. In the percolation vocabulary, $X_\zz=1$ means that site $\zz$ is open (for percolation), and closed otherwise. We denote as $\cI=\{\zz\in\ZZ^2\,:\,X_\zz=1\}$ the set of open sites, and the main question in this oriented percolation model is to characterize the almost sure existence of an infinite path $f\in\cS_R$ such that $f(n)\in\cI$ for all $n\in\bN$. We say that $f\in \cS_R$ is an open crossing of $A_k$ if it is a crossing of $A_k$ and $X_{f(n)}=1$ for all $n\in [0,T_f] \cap \bN$. The multiscale renormalization scheme allow us attack this question by looking at the decay of the probability of the event 
$$D_k\df\left\{\mbox{there exists no open crossing }\,f\in\cS_R\,\mbox{ of }A_k\right\}\,.$$
For $\zz\in\ZZ^2$ let $A_k(\zz)$ and $B_k(\zz)$ be the translated sets of $A_k$ and $B_k$, respectively, and let $D_k(\zz)$ be the event that there exists no open crossing $f\in\cS_R$ of $A_k(\zz)$. It is worth to notice that the event $D_k(\zz)$ has support in the box $B_k(\zz)$, and that its characteristic function is non-increasing (with respect to the configuration of zeros and ones).  

In our context we assume that we can construct simultaneously (coupling) a family of translation invariant random elements $X^\rho:\Omega\to\{0,1\}^{\ZZ^2}$, indexed by $\rho>0$, such that $\cI^\rho\subseteq\cI^{\rho'}$ for $\rho\leq \rho'$, which is equivalent to $X_\zz^\rho\leq X_\zz^{\rho'}$ for all $\zz\in\ZZ^2$ (that induces a partial order in $\{0,1\}^{\ZZ^2}$). We further require a sprinkling decoupling property for this family of random elements, as stated below. The probability law of each random element $X^\rho$ is denoted $\P^\rho$, and expectation $\E^\rho$.    

\begin{hyp}\label{PercolationDecay}
Let $0<a<b$. If $u \in I = (a,b)$, there exist constants $C_I,\epsilon_I>0$ for which the following holds: for all $\epsilon\in(0,\epsilon_I)$ if $f_i:\{0,1\}^{\ZZ^2}\to [0,1] $ for $i=1,2$ are non-increasing functions with support $B_i$ such that 
$$\d(B_1,B_2)\geq C_I\epsilon^{-1}\left(\per(B_1)+\per(B_2)\right)\,,$$
then
$$\E^{\rho}\left[f_1f_2\right]\leq\E^{\rho-\epsilon}\left[f_1\right]\E^{\rho-\epsilon}\left[f_2\right]+H\left(\epsilon,\d(B_1,B_2)\right)\,,$$
where $H$ is non-decreasing in each of the variables with 
\begin{equation}\label{H-decay}
\lim_{x\to\infty} x^7 H(x^{-\delta_0}, x )=0\,,
\end{equation}
for some $\delta_0>0$.
\end{hyp}

We note that, since $H$ is non-decreasing,  if \eqref{H-decay} holds for some $\delta_0>0$ then it also holds for all $\delta\in(0,\delta_0]$. Fix $ \rho>0$ and $\delta>0$. Since $l_k$ has a super-exponential growth, one can define an increasing sequence $(\rho_k)_{k\geq 0}$ with $\rho_0>0$, 
$\rho_{k+1} - \rho_k = c_\rho l_k^{-\delta}$
and $\lim \rho_k = \rho$, where $c_\rho>0$ is a scale factor that only depends on $\rho>0$ and $l_0$. We can now announce the percolation result that we will use to prove Theorem \ref{Detection}.
\begin{thm}\label{Percolation}
Under Hypothesis \ref{PercolationDecay}, let
$$p_k=\P^{\rho_k}\left[D_k\right ]\,.$$ 
There exists $\tilde{k}\geq 1$ such that if (trigger inequality)
$$p_k\leq l_k^{-4}\mbox{ for some }k\geq \tilde{k}\,,$$ 
then 
$$p_n\leq l_n^{-4}\mbox{ for all }n\geq k\,.$$
In particular,
$$\P^{\rho}\left[\mbox{there exists an infinite open path }f\in\cS_R\right]=1\,.$$
\end{thm}

The proof of this theorem is based on the proof of Theorem 2.5 and Theorem 2.9 from \cite{BaTe}, and it starts with a small modification of  Lemma 2.7 of \cite{BaTe}, that relates the events $D_k$ and $D_{k-1}$.
\begin{lem}\label{Scale}
Let $\upsilon\in(0,1/2)$. There exists $k_0\in\bN$ such that, for each $k\geq k_0$, there exists $M_k\subseteq \ZZ^2$ satisfying: 
\begin{enumerate}
    \item $\#M_k\leq 10 l_{k-1}^{1/2}$;
    \item If $D_k$ happens, there exists $\xx,\yy\in M_k$ such that $D_{k-1}(\xx)$ and $D_{k-1}(\yy)$ happen and 
    $$\d\left(B_{k-1}(\xx),B_{k-1}(\yy)\right)\geq \left(\per(B_{k-1}(\xx)+\per(B_{k-1}(\yy)\right)^{3/2-\upsilon}\,.$$
\end{enumerate}

\end{lem}

\noindent{\bf Proof of Lemma \ref{Scale}\,\,}
Following the lines of \cite{BaTe} we define the collection of sets
 $$ \left\{A_{k-1}(\xx_j):0\leq j\leq \frac{L_k+l_k}{L_{k-1}}\right\}\,,$$
with $\xx_j \df j(l_k,L_k)$ and their corresponding reflection with respect to the diagonal of $A_k$. Through concatenation, open crossings of each of these $2\left(\frac{L_k+l_k}{L_{k-1}}+1\right)$ sets produce an open crossing of $A_k$
(see Figure 2.2 of \cite{BaTe}). The same can be said about the collection of sets
$$ \left\{A_{k-1}(\yy_j):0\leq j\leq \frac{L_k}{L_{k-1}}\right\}$$
with $\yy_j = (l_k,L_k) - j(l_{k-1},L_{k-1})$ and their corresponding reflection with respect to the diagonal of $A_k$. As in \cite{BaTe}, we take $M_k$ to be the set of all points $\xx\in \ZZ^2$ such that $A_{k-1}(\xx)$ is in some of the two collections described above. Then (1) follows from Lemma 2.7 in \cite{BaTe}. By computation $(2.24)$ from \cite{BaTe}, we have that 
$$\d(B_{k-1}(\xx),B_{k-1}(\yy))\geq C l_{k-1}^{\frac{3}{2}}k^{-2}\,,$$
for some $C>0$, which implies (2), since $l_k$ has a super-exponential growth and $\per(B_{k-1})\sim l_{k-1}$.

\hfill$\Box$\\

\noindent{\bf Proof of Theorem \ref{Percolation}\,\,}
 Lemma \ref{Scale} allows to bound the probability of $D_{k+1}$ by the probability that at least one element of each of these two previous collection of sets does not admit an open crossing. If we take $\upsilon\in(0,1/2)$ such that $\delta<\min\{1/2-\upsilon,\delta_0\}$,  set $\epsilon_k \df c_\rho l_{k}^{-\delta}$, and use (2) from Lemma \ref{Scale}, then 
 $$\d\left(B_{k}(\xx),B_{k}(\yy)\right)\geq \left(\per(B_{k}(\xx)+\per(B_{k}(\yy)\right)^{3/2-\upsilon}\geq C_I \epsilon_k^{-1}\left(\per(B_{k}(\xx)+\per(B_{k}(\yy)\right)\,,$$
 for large enough $k$. The last inequality gives the condition needed for Hypothesis \ref{PercolationDecay}. Taking $f_1, f_2$ as the non existence of an open crossing for $B_{k}(\xx)$ and $B_{k}(\yy)$ respectively, by (1) 
 from Lemma \ref{Scale} together with the union bound, 
 $$p_{n+1}\leq \# M_{n}^2\left(p_n^2+H(c_\rho l_n^{-\delta},C_I l_n)\right)\,.$$
Thus, by \eqref{H-decay} in Hypothesis \ref{PercolationDecay}, we can take $\tilde{k}\geq k_0$ sufficiently large such that 
 $$100 l_k^{-1} + 100 l_k^7 H(c_\rho l_k^{-\delta},C_I l_k)\leq 1 \,,\,\forall\, k\geq \tilde{k}\,.$$
 By induction, if $p_{k}l_k^{4}\leq 1$ for some $k\geq \tilde{k}$, then for $n\geq k$ we have that 
 \begin{align*}
     p_{n+1} l_{n+1}^4  &\leq \left(10 l_{n-1}^{1/2}\right)^2 l_{n+1}^4 \left( l_n^{-8} +H(c_\rho l_n^{-\delta},C_I l_n) \right)\\
     &\leq 100 l_n \left( l_n^{3/2}\right)^4\left( l_n^{-8} +H(c_\rho l_n^{-\delta},C_I l_n) \right)\\
     &\leq 100 l_n^{-1} + 100 l_n^7 H(c_\rho l_n^{-\delta},C_I l_n)\\
     &\leq  1.
 \end{align*}
Finally, as in \cite{BaTe}, we notice that by Borel-Cantelli there is an open crossing of $A_k(\xx)$ for all $\xx \in M_k$ for all sufficiently large $k$. The existence of an open crossing of $A_k(\xx)$ at $\rho_k$ implies the existence at $\rho$, as $\rho_k < \rho$. By concatenation this implies the existence of an infinite open path $f\in\cS_R$ for the parameter $\rho$.
 
\hfill$\Box$\\

\subsection{Proof of Theorem \ref{Detection}}

For the sake of simplicity, we take $\Delta_s=\Delta_t=1$. For $r>0$ let $I_{x,r}\df (x-r,x+r)$. Given $\rho>0$ let $\lambda\df \rho^{-1/2}>0$ (so that $\rho=\lambda^{-2}$). To link the detection problem to a percolation model we define for $(x,n)\in\bZ^2$ 
$$X^\rho_{(x,n)}=\begin{cases}  
1 \mbox{ (open) }\,,\mbox{ if }\eta^\lambda_t(I_{x,r})=0\,\,\forall\,t\in[n,n+1)\,,\\
0 \mbox{ (closed) }\,,\,\mbox{ otherwise }\,.
\end{cases} $$ 
If $X^\rho_{(x,n)}=1$ then a target at $x$ is not detected by Hammersley's process for the period of time $t\in[n,n+1)$. Recall that $\cI^\rho$ denotes the collection of open sites, and that, and by attractiveness of Hammersley's process, we can take $(\cI^\rho)_{\rho>0}$ as an increasing family of subsets of $\bZ^2$. Then percolation of the set $\cI^\rho$ using the set of paths $\cS_R$ implies non-detection with positive probability. We have to verify that this percolation system  satisfies the necessary conditions. By space-time stationarity of Hammersley's process, it is clearly translation invariant, and Hypothesis \ref{PercolationDecay} follows from Theorem \ref{Sprinkling}. To obtain the trigger inequality we set again $\lambda_k\df \rho_k^{-1/2}$ (so that $\rho_k=\lambda_k^{-2}$), and emphasize the following geometrical remarks:
\begin{itemize}
    \item[(i)] For $N_k=L_k+\lfloor l_k\rfloor +1$ (which bounds the size of the jumps), to obtain an open $\cS_{N_k}$-crossing of $A_k$ it is sufficient and necessary to have an open vertical $\cS_{N_k}$-crossing of $[0,l_k]\times[0,L_k]$ since, in such a case, one can do one step of size $N_k$ to exit on the right side of the  top part of $A_k$;
    \item[(ii)] If at every horizontal level $[0,l_k]\times \{j\}$ with $j\in[0,L_k)$ we have a open site, than one can build an open vertical $\cS_{N_k}$-crossing of $[0,l_k]\times[0,L_k]$ by jumping from level to level over open sites (recall that $N_k>l_k$ is large enough for that).      
\end{itemize}
Now, if $\eta^{\lambda_k}_n(I_{x,r})=0$ and there is no $\omega$ space-time Poisson clock in  $I_{x,r}\times [n,n+1)$ then $X^{\rho_k}_{(x,n)}=1$. Hence 
$$\{X^{\rho_k}_{(x,n)}=0\}\subseteq J_{(x,n)}\df\{\eta^{\lambda_k}_n(I_{x,r})>0\mbox{ or }\omega\left(I_{x,r}\times [n,n+1)\right)>0\}\,.$$
A simple calculation shows that 
$$\P^{\lambda_k}\left[J_{(x,n)}\right]= 1- e^{-2r\lambda_k}e^{-2r} =1-e^{-2r(\lambda_k+1)}\leq 1-e^{-2r(\lambda_0+1)}\,,$$
since $\lambda_k\leq \lambda_0$ (recall that $0<\rho_0\leq \rho_k$). If $|x-y|\geq 2r$ then $I_{x,r}\cap I_{y,r}=\emptyset$ which implies that $J_{(x,n)}$ and $J_{(y,n)}$ are independent. Choosing $N_k=L_k+\lfloor l_k\rfloor+1$, and using remarks (i) and (ii), 
\begin{multline*}
p_k=\P^{\rho_k}\left[\mbox{ no open vertical $\cS_{N_k}$-crossing of $[0,l_k]\times[0,L_k]$ }\right]\\
\leq (L_k+1)\P^{\lambda_k}\left[\bigcap_{x\in[0,l_k]\cap 2r\bN}J_{(x,0)}\right]\leq (L_k+1)\left( 1-e^{-2r(\lambda_0+1)}\right)^{\lfloor l_k/2r\rfloor}\\
\leq (2l_k+1)\left( 1-e^{-2r(\lambda_0+1)}\right)^{\lfloor l_k/2r\rfloor}\,.
\end{multline*}
This inequality allows us to choose $\tilde k\geq 1$ such that $l^4_{\tilde k}p_{\tilde k}\leq 1$, which implies the trigger inequality. 

\section{Multiscale Analysis for the Radom Walk}

For the random walk analyses we parallel the proof of Proposition 3.6 of \cite{HiKiTe}, and set different collections of scales $\{L_k\,:\,k\geq k\}$ and parameters $\{\rho_k\,:\,k\geq 1\}$. Define recursively,
$$L_0\df 10^{10}\mbox{ and }L_{k+1}\df l_k L_k\,,\,\forall \,k\geq 0\,,$$
where $l_k\df \lfloor L_k^{1/4}\rfloor$. For $L\geq 1$ and $h\geq 1$, define 
$$B_{L}^h\df [-hL,2hL)\times [0,hL)\subseteq \R^2 \,$$
and $B_{L}^h(\w)\df \w+B_{L}^h$. Let $0<\rho<\rho_{c-}$ and choose $\delta>0$ such that $\rho<\rho_{c-}-2\delta$. Let $\epsilon_k\df L_k^{-1/16}$, fix $k_6$ so that $\sum_{k\geq k_6}\epsilon_k\leq \frac{\delta}{4}$, 
and define the sequence of parameters 
$$\rho_{k_6}\df \rho_{c-} -\frac{5}{4}\delta\,,\mbox{ and }\rho_{k+1}\df \rho_k-\epsilon_k\,,\,\forall\,k\geq k_6\,.$$
Thus 
$$0<\rho<\rho_{c-}-2\delta< \rho_\infty\df\lim_{k\to\infty}\rho_k\,.$$ 
By Theorem \ref{Sprinkling}, if $f_1$ and $f_2$ are non-increasing functions taking values in $[0,1]$, supported on $B_{L_k}^h(\w_1)$ and  $B_{L_k}^h(\w_2)$, respectively, such that  
$$|\pi_2(\w_2)-\pi_2(\w_1)|\geq C_1\epsilon_k^{-1} h L_k\,$$
then 
\begin{equation}\label{MultiScaleSprink}
\E^{\rho_{k+1}}[f_1f_2]\leq \E^{\rho_{k}}[f_1]\E^{\rho_{k}}[f_2]+ C_2  e^{-c_1 (hL_k)^{?}}\,.
\end{equation}

Recall that $l_k = \lfloor L_k^{1/4}\rfloor$. Since $\rho<\rho_{c-}-2\delta<\rho_{c-}-\delta<v_-(\rho_{c-}-\delta)$, we can choose $k_7\geq k_6$ large enough such that 
$$\sum_{k\geq k_7}\frac{ 2C_1 L^{1/16}_k}{l_k}<\frac{v_-(\rho_{c-}-\delta)-\rho}{4}\,.$$
Define a decreasing sequence of velocities by setting 
$$\tilde v_{k_7}\df v_-(\rho_{c-}-\delta)-\frac{1}{8}\left(v_-(\rho_{c-}-\delta)-\rho\right)\,,$$
and
$$ \tilde v_{k+1}\df \tilde v_k-\frac{ 2C_1 L^{1/16}_k}{l_k}\,\mbox{ for }k\geq k_7\,.$$
Thus,
$$\tilde v_\infty\df \lim_{k\to\infty}\tilde v_k\geq 
\rho+\frac{5}{8}\left(v_-(\rho_{c-}-\delta)-\rho\right)\,,$$
and $\tilde{v}_k\in(0,1)$ for all $k\geq k_7$. The next lemma is a simple adaptation of Lemma 7.2 of \cite{HiKiTe}, where now the distance between the support of the events must be slightly larger than the symmetric case, in order to meet the requirements of the sprinkling decoupling inequality \eqref{MultiScaleSprink}. 

\begin{lem}\label{SprinkInduction}
Let $h\geq 1$, $k\geq k_7$ and $m\in M_{k+1}^h$. If $\tilde{A}_m(\tilde{v}_{k+1})$ holds, then there exist $m_1=(h,k,\w_1),m_2=(h,k,\w_2)\in C_m$ such that $\pi_2(\w_2)-\pi_2(\w_1)\geq C_1\epsilon_k^{-1} h L_k$ and both $\tilde{A}_{m_1}(\tilde{v}_{k})$ and $\tilde{A}_{m_2}(\tilde{v}_{k})$ hold.    
\end{lem}

\noindent{\bf Proof of Lemma \ref{SprinkInduction}\,\,} 
Let $m\in M_{k+1}^h$ and assume that $\tilde{A}_m(\tilde{v}_{k+1})$ holds. Then there exists $\y\in[0,hL_{k+1})\times \{0\}$ such that \begin{equation}\label{Induc1}
X^\y_{hL_{k+1}}-\pi_1(\y)\leq \tilde{v}_{k+1}hL_{k+1}\,.
\end{equation}
For any $i\in\{0,\dots,l_k\}$ there exists a unique $m'=m'(i)\in C_m$ such that $X^{\y}_{ihL_k}\in I_{m'}$, and let $M(X^\y)$ denote the collection of these random indices $m'(i)$ for $i\in\{0,\dots,l_k\}$. If we have at least  $C_1\epsilon_k^{-1}+1$ indices $m'\in M(X^\y)$ such that $A_{m'}(\tilde{v}_k)$ holds, then by picking $m_1$ and $m_2$ the first and the last indices (along the path), respectively, then we certainly have that $\pi_2(\w_2)-\pi_2(\w_1)\geq C_1\epsilon_k^{-1} h L_k$.     

By contradiction, assume that we have at most $C_1\epsilon_k^{-1}$ indices $m'\in M(X^\y)$ such that $\tilde A_{m'}(\tilde{v}_k)$ occurs. Now, for each $i\in\{0,\dots,l_k\}$ and $m'=m'(i)$ we have $X^\y_{ih L_k}\in I_{m'}$ and if  $\tilde A_{m'}(\tilde{v}_k)$ does not occur then  
$$X^\y_{(i+1)hL_k}-X^\y_{ihL_k}\geq \tilde{v}_k h L_k\,. $$
Otherwise, we can bound the increment using the deterministic bound
$$X^\y_{(i+1)hL_k}-X^\y_{ihL_k}\geq -hL_k\,. $$
Therefore,
\begin{eqnarray*}
X^\y_{hL_{k+1}}-\pi_1(\y)&=&\sum_{i=0}^{l_k}X^\y_{(i+1)hL_k}-X^\y_{ihL_k}  \\
&\geq& - C_1\epsilon_k^{-1} hL_k +(l_k-C_1\epsilon_k^{-1})\tilde{v}_kh L_k\\
&=&hL_{k+1}\left( \tilde{v}_k -\frac{C_1\epsilon_k^{-1}(1+\tilde{v}_k)}{l_k} \right)\,.
\end{eqnarray*}
In order to have that 
$$ \tilde{v}_k -\frac{C_1\epsilon_k^{-1}(1+\tilde{v}_k)}{l_k}=\tilde{v}_k -\frac{C_1 L_k^{1/16}(1+\tilde{v}_k)}{l_k}>\tilde{v}_{k+1}=\tilde{v}_k- \frac{ 2C_1 L^{1/16}_k}{l_k}\,,$$
we only need that $C_1(1+\tilde{v}_k)<2C_1$, which is  true since $v_k\in (0,1)$ for all $k\geq k_7$. This yields to 
$$ X^\y_{hL_{k+1}}-\pi_1(\y)> \tilde{v}_{k+1} h L_{k+1}\,,$$
which contradicts \eqref{Induc1}.

\hfill$\Box$\\ 

The events $\tilde{A}_m(v)$ are non-increasing with respect to the environment (and so non-decreasing in $\rho$). By combining Lemma \ref{SprinkInduction} together with \eqref{MultiScaleSprink} and the union, we have that (induction step)
\begin{equation}\label{Induc2}
\tilde{p}_{hL_{k+1}}(\tilde{v}_{k+1},\rho_{k+1})\leq 9 l_k^4\left(\tilde{p}_{hL_{k}}(\tilde{v}_{k},\rho_{k})^2+C_2e^{-c_1(hL_k)^?}\right)\,,\,\forall\,k\geq k_7\,.    
\end{equation}
Choose $k_8\geq k_7$ such that
$$9l_k^4\left( e^{-4(2-(5/4)^{3/2}(\log L_k)^{3/2})}+C_2e^{-c_1 L_k^?+8(\log L_k)^{3/2} } \right)\leq 1\,,\,\forall\,k\geq k_8\,.$$
Since $\tilde{v}_{k_8}<v_-(\rho_{c-}-\delta)$ and $\rho_{k_8}<\rho_{c-}-\delta$, we have that  
$$\liminf_{H\to\infty}\tilde{p}_H(\tilde{v}_{k_8},\rho_{k_8})\leq \liminf_{H\to\infty}\tilde{p}_H(\tilde{v}_{k_8},\rho_{c-}-\delta)=0\,.$$
Therefore, we can choose $h_0\geq 1$ large enough so that (trigger step) 
$$\tilde{p}_{h_0 L_{k_8}}(\tilde{v}_{k_8},\rho_{k_8})\leq e^{-4(\log L_{k_8})^{3/2}} \,.$$
By \eqref{Induc2}, if we assume that 
$$\tilde{p}_{h_0 L_{k}}(\tilde{v}_{k},\rho_{k})\leq e^{-4(\log L_{k})^{3/2}}\,,$$
for $k\geq k_8\geq k_7$ then 
\begin{eqnarray*}
\frac{\tilde{p}_{hL_{k+1}}(\tilde{v}_{k+1},\rho_{k+1})}{e^{-4(\log L_{k+1})^{3/2}}}&\leq & 9 l_k^4e^{4(\log L_{k+1})^{3/2}}\left(\tilde{p}_{hL_{k}}(\tilde{v}_{k},\rho_{k})^2+C_2e^{-c_1(hL_k)^?}\right) \\
&\leq &9l_k^4\left( e^{-4(2-(5/4)^{3/2}(\log L_k)^{3/2})}+C_2e^{-c_1 L_k^?+8(\log L_k)^{3/2} } \right)\\
&\leq & 1\,.
\end{eqnarray*}
Hence, by induction,
$$\tilde{p}_{h_0 L_{k}}(\tilde{v}_{k},\rho_{k})\leq e^{-4(\log L_{k})^{3/2}}\,,\,\forall\, k\geq k_8\,.$$
Since $\tilde{A}_m(v)$ is non-decreasing in $v$ and $\rho$, $\rho + \frac{5}{8}\left(v(\rho_{c-}-\delta)-\rho\right)<\tilde{v}_k$ and $\rho<\rho_{c-}-2\delta<\rho_k$, we have that 
\begin{equation}\label{MultiScaleIneq}
\tilde{p}_{h_0 L_{k}}\left( \rho + \frac{5}{8}\left(v(\rho_{c-}-\delta)-\rho\right) ,\rho\right)\leq e^{-4(\log L_{k})^{3/2}}\,,\,\forall\, k\geq k_8\,. 
\end{equation}

From here we have to check the interpolation argument to conclude the proof of the first part of Proposition \ref{Ballistic} (see the proof of Proposition 3.6 in \cite{HiKiTe}). For large enough $H \geq 1$ let us define $\bar{k}\geq k_8$ as the integer that satisfies:
\begin{equation*}
    h_0 L_{\bar{k}+1} \leq H \leq  h_0 L_{\bar{k}+2}
\end{equation*}
Therefore,
\begin{equation*}
    \tilde{p}_{h_0 L_{\bar{k}}}\left(\rho + \frac{5}{8}\left(v(\rho_{c-}-\delta)-\rho\right) ,\rho\right)\leq e^{-4(\log L_{\bar{k}})^{3/2}}\,.
\end{equation*}

Fix some \( w \in \mathbb{R}^2 \) pave the box \( B_H(w) \) with boxes \( B_m \), where \( m \in M_{h_0 \bar{k}} \) is given by:
\[
m = (h_0, \bar{k}, w + (x h_0 L_{\bar{k}}, y h_0 L_{\bar{k}})),
\]
where the indices \( x \) and \( y \) satisfy:
\[
-\left\lceil \frac{H}{h_0 L_{\bar{k}}} \right\rceil \leq x \leq \left\lceil \frac{2H}{h_0 L_{\bar{k}}} \right\rceil, \quad 0 \leq y \leq \left\lceil \frac{H}{h_0 L_{\bar{k}}} \right\rceil.
\]
Let \( M \) be the set of such indices. Note that:
\[
|M| \leq 4 \left( \frac{H}{h_0 L_{\bar{k}}} \right)^2 \leq 4 \left( \frac{L_{\bar{k}+2}}{L_{\bar{k}}} \right)^2 \leq 4 \left(\frac{L_{\bar{k}}^{25/16}}{L_{\bar{k}}}\right)^2  \leq 4 L_{\bar{k}}^{9/8}.
\]
Let $v_0 = \rho + \frac{5}{8}\left(v(\rho_{c-}-\delta)-\rho\right)$. On the event \( \bigcap_{m \in M} (\tilde{A}_m(v_0))^c \), for any \( \y \in I_H^1(w) \), the displacement of \( X^\y \) up to time $H$ can be bounded as:
\begin{align*}
    X^\y (H) -  \pi_1(y)  &=  X^\y (H) - X^\y (H / h_0 L_{\bar{k}} \rfloor h_0 L_{\bar{k}}) + X^\y(\lfloor H / h_0 L_{\bar{k}} \rfloor h_0 L_{\bar{k}}) - \pi_1(y)\\
    &\geq     -h_0 L_{\bar{k}} + v_0 \lfloor H / h_0 L_{\bar{k}} \rfloor h_0 L_{\bar{k}} \\
    &\geq     -(1+v_0)h_0 L_{\bar{k}} + v_0 H,
\end{align*}
where we used the fact that \( \tilde{A}_m(v_0) \) does not occur for any \( m \in M \). Therefore, for any \( y \in I_H^1(w) \), we have:
\[
\frac{X^\y(H) - \pi_1(y)}{H} \geq -\frac{(1+v_0)h_0 L_{\bar{k}}}{H} + v_0 \geq \tilde{v},
\]
for any fixed $\rho<\tilde{v} < v_0$ by taking $H$ sufficiently large.
Thus,
\[
\tilde{p}_h(\tilde{v},\rho) \leq 4L_{\bar{k}}^{9/8} \exp\left(-4 (\ln L_{\bar{k}})^{3/2} \right) \leq e^{-2 (\ln H)^{3/2}}.
\]
The conclusion of Theorem 3 is then obtain by taking the supremum over all $ \omega \in [0, 1) \times \{0\}$ and then properly choosing the constant $\tilde{c}_1$ in order to accommodate small values of $H$.

\appendix
\section{Tail Bounds for Poisson Random Variables}\label{Poisson} 
Recall that $\Poi^{\lambda}$ denotes a Poisson random variable with parameter $\lambda$. The proof of Lemma \ref{PoissonConcentration} below was based on Cl\'ement Canonne's notes on Poisson tail bounds (see  \href{https://github.com/ccanonne/probabilitydistributiontoolbox/blob/master/poissonconcentration.pdf}{GitHub}), which in turn were derived from David Pollard's unpublished notes \cite{Po}. For the sake of completeness, we decide to include the proof in this article.   
\begin{figure}
    \centering
    \subfloat{{\includegraphics[width=7cm]{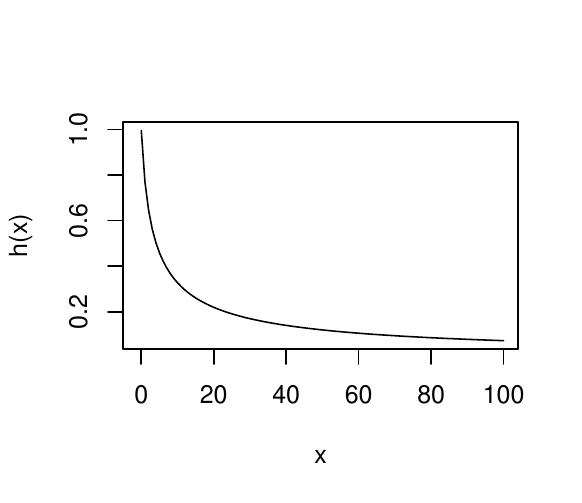} }}
    \qquad
    \subfloat{{\includegraphics[width=7cm]{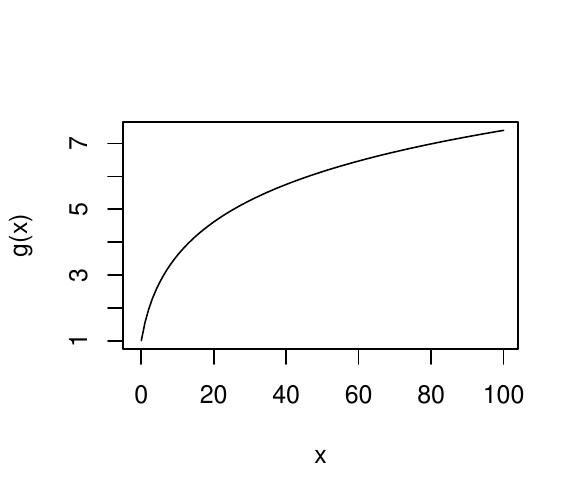} }}
    \caption{A plot using R of the functions that appear in the proof of Lemma \ref{PoissonConcentration}.}%
    \label{FigMGFRPlot}
\end{figure}

\begin{lem}\label{PoissonConcentration}
We have that  
$$\max\big\{\P\left[\Poi^{\lambda}\geq \lambda +x\right],\P\left[\Poi^{\lambda}\leq \lambda -x\right]\big\}\leq e^{-\frac{x^2}{2(x+\lambda)}}\,,\,\forall\,x>0\,.$$ 
\end{lem}
\noindent{\bf Proof of Lemma \ref{PoissonConcentration}\,\,} 
The moment generating function method applied to the Poisson distribution \cite{Po} implies that 
\begin{equation}\label{PoissonLD1}
\P\left[\Poi^{\lambda}\geq \lambda +x\right]\leq e^{\inf_{\theta>0}\{\lambda(e^\theta -1)-\theta(\lambda+x)\}}=e^{-\frac{x^2}{2\lambda}h\left(\frac{x}{\lambda}\right)}\,\,\mbox{  for }x>0\,,
\end{equation}
and 
\begin{equation}\label{PoissonLD2}
\P\left[\Poi^{\lambda}\leq \lambda -x\right]\leq e^{\inf_{\theta>0}\{\lambda(e^{-\theta} -1)+\theta(\lambda-x)\}}= e^{-\frac{x^2}{2\lambda}h\left(-\frac{x}{\lambda}\right)}\,\,\mbox{  for }0<x<\lambda\,,
\end{equation}
where $h:[-1,\infty)\to \R$ is the function defined as 
$$h(x)\,=\,\begin{cases} 2\frac{(1+x)\log(1+x)-x}{x^2}\,\,\mbox{ for }x\neq -1,0\\
2\quad\quad\quad\quad\quad\quad\,\,\,\,\,\mbox{ for }x= -1\\
1\quad\quad\quad\quad\quad\quad\,\,\,\,\,\mbox{ for }x= 0\end{cases}\,=\,\int_0^1\frac{2(1-s)}{1+sx}ds\,.$$ 
To justify the second equality above \cite{Po}, one can use the integral form of Taylor's theorem,  
$$f(x)-f(0)-xf'(0)=\int_0^x(x-t)f''(t)dt=\frac{1}{2}x^2\int_0^1 2(1-s)f''(sx)ds\,,$$
with $f(x)=(1+x)\log(1+x)$. Hence, $h$ is continuous, decreasing, non-negative, and $g(x)=(1+x)h(x)$ is increasing (Figure \ref{FigMGFRPlot}), which can be deduced straightforward by using the integral form of $h$. Therefore, $h(x)\geq h(0)/(1+x)=1/(1+x)$ for $x\geq 0$, and   
$$\frac{x^2}{2\lambda}h\left(\frac{x}{\lambda}\right)\geq \frac{x^2}{2(x+\lambda)}\,\implies\,\P\left[\Poi^{\lambda}\geq \lambda +x\right]\leq e^{-\frac{x^2}{2\lambda}h\left(\frac{x}{\lambda}\right)}\leq e^{-\frac{x^2}{2(x+\lambda)}}\,.$$
For the other bound, we have that 
$$\P\left[\Poi^{\lambda}\leq \lambda -x\right]\leq e^{-\frac{x^2}{2\lambda}h\left(-\frac{x}{\lambda}\right)}\leq e^{-\frac{x^2}{2\lambda}h(0)}\leq e^{-\frac{x^2}{2(\lambda+x)}}\,,$$
for $0<x<\lambda$. Notice that for $x \geq\lambda$ the upper bound is trivial. 

\hfill$\Box$\\

\end{document}